\documentclass[11pt,a4paper]{article}
\usepackage{graphicx}%pictures
\usepackage{amsmath}
\usepackage{amssymb}
\usepackage{amsthm}
\usepackage{geometry}
\usepackage{fancyhdr}%Heads
\usepackage{color} %colors
\usepackage{overpic}%pictures with text
\usepackage{mathabx}

\usepackage[all]{xy}%commutative graph
\newcommand{\al}{\alpha}

\newcommand{\be}{\beta}
\newcommand{\ga}{\gamma}

\newcommand{\ra}{\rightarrow}%arrows

\newcommand{\xra}{\xrightarrow}

\newcommand{\bpf}{\begin{proof}}%thm,prop,lem,pf,example,remark
\newcommand{\epf}{\end{proof}}
\newcommand{\bthm}{\begin{thm}}
\newcommand{\ethm}{\end{thm}}
\newcommand{\bprop}{\begin{prop}}
\newcommand{\eprop}{\end{prop}}
\newcommand{\bcor}{\begin{cor}}
\newcommand{\ecor}{\end{cor}}
\newcommand{\blem}{\begin{lem}}
\newcommand{\elem}{\end{lem}}
\newcommand{\bdefn}{\begin{defn}}
\newcommand{\edefn}{\end{defn}}
\newcommand{\bexmp}{\begin{exmp}}
\newcommand{\eexmp}{\end{exmp}}
\newcommand{\brem}{\begin{rem}}
\newcommand{\erem}{\end{rem}}

\newcommand{\bdia}{\begin{displaymath}\xymatrix}
\newcommand{\edia}{\end{displaymath}}
\newcommand{\beq}{\begin{equation*}\begin{aligned}}
\newcommand{\eeq}{\end{aligned}\end{equation*}}

\newcommand{\intg}{\mathbb{Z}}%numbers

\newtheorem{thm}{\textbf {Theorem}}[section]
\newtheorem{cor}[thm]{\textbf{Corollary}}
\newtheorem{prop}[thm]{\textbf{Proposition}}
\newtheorem{lem}[thm]{\textbf{Lemma}}
\newtheorem{conj}[thm]{Conjecture}

\newtheorem{quest}[thm]{Question}

\theoremstyle{definition}
\newtheorem{defn}[thm]{\textbf{Definition}}

\newtheorem{exmp}[thm]{Example}

\theoremstyle{remark}
\newtheorem{rem}[thm]{Remark}

\title{Contact structures, excisions, and sutured monopole Floer homology}

\author{Zhenkun Li}
\date{}

\begin{document}
\bibliographystyle{plain}%for reference

\maketitle
\begin{abstract}
In this paper, we explore the interplay between contact structures and sutured monopole Floer homology. First, we study the behavior of contact elements, which were defined by Baldwin and Sivek in \cite{baldwin2016contact}, under the operation of performing Floer excisions, which was introduced to the context of sutured monopole Floer homology by Kronheimer and Mrowka in \cite{kronheimer2010knots}. We then compute the sutured monopole Floer homology of some special balanced sutured manifolds, using tools closely related to contact geometry. For application, we obtain an exact triangle for the oriented skein relation in monopole theory and derive a connected sum formula for sutured monopole Floer homology.
\end{abstract}

\tableofcontents%table of contents
\newpage

%————Start from here————
\section{Introduction}
Sutured monopole and instanton Floer homologies were introduced by Kronheimer and Mrowka in \cite{kronheimer2010knots}. They were designed to be the counterparts of Juh\'asz's sutured (Heegaard) Floer homology \cite{juhasz2006holomorphic} in the monopole and the instanton settings, respectively. 

It was shown by Kutluhan, Lee, and Taubes in \cite{kutluhan2010hf} and in subsequent papers, that monopole Floer homology is isomorphic to Heegaard Floer homology. Using their work, in \cite{baldwin2016equivalence}, Baldwin and Sivek \cite{baldwin2016equivalence} (or Lekili in \cite{lekili2013heegaard}) proved that sutured monopole Floer homology and sutured (Heegaard) Floer homology are isomorphic to each other. Therefore, if we simply aim at computing monopole Floer homologies, we could make use of the isomorphism and look at the Heegaard Floer side, which is known for being more computable. However, there are no known isomorphisms between instanton Floer homology and any other version of Floer homology. So, it is still valuable to develop some techniques to compute (sutured) monopole Floer homology, which could also shed light on computing (sutured) instanton Floer homology. In this paper, we will stay in the monopole realm, and it is worth mentioning that most of the time we only make use of the formal properties of monopole Floer homology. Thus, the same argument can easily be adapted to instanton theory and potentially to any other Floer theory that shares the same set of formal properties.

A balanced sutured manifold is a compact oriented $3$-manifold $M$ whose boundary is divided by a (possibly disconnected) closed oriented curve $\ga$, which is called the suture, into two parts of the same Euler characteristics. To define monopole Floer homology on such a manifold, one needs to construct a closed oriented $3$-manifold $Y$, together with a closed oriented surface $R\subset Y$, out of $(M,\ga)$. To achieve this, one needs to glue a thickened surface $T\times[-1,1]$ to $M$ along an annular neighborhood of the suture and to identify the two boundary components of the resulting $3$-manifold. Here, $T$ is called an auxiliary surface and is chosen so that it is connected and oriented. We also require that its boundary is identified with the suture. The pair $(Y,R)$ is called a closure. Also, choose a non-separating curve $\eta\subset R$ to construct local coefficients. Then, define the sutured monopole Floer homology of the balanced sutured manifold $(M,\ga)$ to be:
$${SHM}(M,\ga):=HM(Y|R;\Gamma_{\eta}):=\bigoplus_{c_1(\mathfrak{s})[R]=2g(R)-2}\widecheck{HM}(Y,\mathfrak{s};\Gamma_{\eta}).$$

If $(M,\ga)$ is equipped with a contact structure $\xi$ so that $\partial{M}$ is convex and $\ga$ is the dividing set, then, in \cite{baldwin2016contact}, Baldwin and Sivek introduced a way to extend the contact structure $\xi$, which is on $M$, to a contact structure $\bar{\xi}$ on $Y$. Here, $Y$ is a suitable closure of $(M,\ga)$. Hence, by Kronheimer, Mrowka, Ozsv\'ath, and Szab\'o \cite{kronheimer2007monopolesandlens}, one can define the contact invariant 
$$\phi_{\xi}=\phi_{\bar{\xi}}\in HM(-Y|-R;\Gamma_{-\eta})={SHM}(-M,-\ga)$$
in sutured monopole Floer theory. 

Contact structures and contact elements have played very important roles in sutured (Heegaard) Floer theory. The constructions of gluing maps and cobordism maps both need contact structures (see \cite{juhasz2016cobordisms,honda2008contact}). The re-construction of $HFK^{-}$ using direct systems of sutured Floer homologies by Etnyre, Vela-Vick, and Zarev in \cite{etnyre2017sutured} also involves contact geometry in an essential way. Moreover, in \cite{kalman2017tight}, K\'alm\'an and Mathews provided some examples where the generators of the sutured (Heegaard) Floer homologies of some families of balanced sutured manifolds are in one-to-one correspondence to the tight contact structures on those manifolds.

In this paper, we will go deeper into the interplay between contact structures and sutured monopole Floer theory. We have two main topics.

\subsection{Contact elements and Floer excisions}
First, we study the behavior of contact elements under Floer excisions. In \cite{kronheimer2010knots}, Kronheimer and Mrowka used both connected and disconnected auxiliary surfaces to define sutured monopole Floer homology, as well as proving many basic properties of it. The isomorphism between sutured monopole Floer homologies arising from connected and disconnected auxiliary surfaces was constructed through Floer excisions. Later, Baldwin and Sivek constructed the contact invariants by only using connected auxiliary surfaces. So, it is interesting to ask whether their construction can be extended to the case of using disconnected auxiliary surfaces, and how contact elements arising from connected and disconnected auxiliary surfaces are related by Floer excisions. The answers to these questions can help us to better understand the maps associated to the trace and co-trace cobordisms, as well as the behavior of contact elements under proper sutured manifold decompositions.

Suppose that for $i=1,2$, $(M_i,\ga_i)$ is a balanced sutured manifold, and $T_i$ is a connected auxiliary surface, which gives rise to a closure $(Y_i,R_i)$ of $(M_i,\ga_i)$. If we cut $T_1$ and $T_2$ along one suitably chosen non-separating simple closed curve on each of them and glue the resulting surfaces together along the newly created boundary components, then we obtain a connected surface $T$. We can use $T$ as an auxiliary surface to close up $(M_1\sqcup M_2,\ga_1\cup\ga_2)$ and obtain a connected closure $(Y,R)$. In \cite{kronheimer2010instanton}, Kronheimer and Mrowka constructed a Floer excision cobordism $W$ from $(Y_1\sqcup Y_2)$ to $Y$, and after choosing suitable curves on $R_1$, $R_2$ and $R$ to support local coefficients, the cobordism $W$ will induce a map
\begin{equation}\label{eq_excision_map}
F=HM(-W):HM(-(Y_1\sqcup Y_2)|-(R_1\cup R_2);\Gamma_{-(\eta_1\cup \eta_2)})\ra HM(-Y|-R;\Gamma_{-\eta}).
\end{equation}

Suppose further that for $i=1,2$, $(M_i,\ga_i)$ is equipped with a contact structure $\xi_i$ so that $\partial{M}_i$ is convex and $\ga_i$ is the dividing set. Then, by Baldwin and Sivek \cite{baldwin2016contact}, there are corresponding contact structures $\bar{\xi}_1$, $\bar{\xi}_2$ and $\bar{\xi}$ on $Y_1$, $Y_2$ and $Y$ respectively. Then, there are contact elements $\phi_{\bar{\xi}_1\cup\bar{\xi}_2}$ and $\phi_{\bar{\xi}}$ associated to them. In this paper, we prove the following: 
\begin{thm}\label{thm_1}
The map $F$ in (\ref{eq_excision_map}) preserves contact elements up to multiplication by a unit:
$$F(\phi_{\bar{\xi}_1\cup\bar{\xi}_2})\doteq \phi_{\bar{\xi}}.$$
Here $\doteq$ means equal up to multiplication by a unit.
\end{thm}

\brem
As explained in Baldwin and Sivek \cite{baldwin2016contact} and Krongeimer and Mrowka \cite{kronheimer2010knots}, to define contact elements and to carry out Floer excisions along tori, it is necessary to use local coefficients. Throughout this paper, we always use the Novikov Ring or the mod 2 Novikov ring as the base ring to construct local coefficients.
\erem

The result in Theorem \ref{thm_1}, however, is not fully satisfactory. Suppose that $(M_1,\ga_1)$ and $(M_2,\ga_2)$ are the same as in Theorem \ref{thm_1}, and $(M,\ga)$ is a connected sutured manifold so that there is a diffeomorphism
$$g:\partial{M}\xra{\cong}\partial{M}_1\sqcup\partial{M}_2,$$
which sends $\ga$ to $\ga_1\cup \ga_2$. Then, we can use either $T_1\sqcup T_2$ or $T$ to construct a closure of $(M,\ga)$. The two resulting closures are still related by a Floer excision cobordism $W$, and $W$ induces a map between the corresponding monopole Floer homologies, which is similar to the one in (\ref{eq_excision_map}). The proof of Theorem \ref{thm_1} in this paper, however, fails to cover the case when $(M,\ga)$ is connected. We make the following conjecture:

\begin{conj}\label{conj_2}
Theorem    \ref{thm_1} continues to hold when we replace the disjoint union $(M_1\sqcup M_2,\ga_1\cup\ga_2)$ by a connected balanced sutured manifold $(M,\ga)$, as described in the above parograph.
\end{conj}

Some evidence or ideas of the proof can be found in \cite{niederkruger2011weak}. In \cite{niederkruger2011weak}, Niederkruger and Wendl introduced an operation called slicing, which coincides with performing a Floer excision, and an operation of attaching torus $1$-handles, which coincides with constructing a Floer excision cobordism $W$. Thus, the cobordism $W$ is equipped with a weak symplectic structure, as explained in \cite{niederkruger2011weak}. Compared with the known results reached by Hutchings and Taubes in \cite{hutchings2013proof} and by Echeverria in \cite{echeverria2018naturality} that both exact symplectic and strong symplectic cobordisms preserve contact elements, we make the following conjecture.

\begin{conj}\label{conj_3}
Suppose $(W,\omega)$ is a weakly symplectic cobordism from $(Y_1,\xi_1)$ and $(Y_2,\xi_2)$. Then, suppose that for $i=1,2$, there is a $1$-cycles $\eta_i\subset Y_i$, so that $\eta_i$ is dual to $\omega|_{Y_i}$. Continuing,  suppose $\nu\subset W$ is a $2$-cycle so that $\partial{\nu}=-\eta_1\cup\eta_2$. Then, the map
$$\widecheck{HM}(W,\mathfrak{s_{\omega}};\Gamma_{\nu}):\widecheck{HM}(-Y_2,\mathfrak{s}_{\xi_2};\Gamma_{-\eta_2})\ra \widecheck{HM}(-Y_1,\mathfrak{s}_{\xi_1};\Gamma_{-\eta_1})$$
will preserve the contact elements:
$$\widecheck{HM}(W,\mathfrak{s_{\omega}};\Gamma_{\nu})(\phi_{\xi_2})\doteq\phi_{\xi_1}.$$
\end{conj}

The confirmation of Conjecture \ref{conj_3} would possibly provide a proof of Conjecture \ref{conj_2}.

\subsection{Connected sum formula}
 In the second half of the paper, we prove the connected sum formula for sutured monopole Floer homology. In particular, we prove the following theorem. 
\bthm
Suppose $(M_1,\ga_1)$ and $(M_2,\ga_2)$ are two balanced sutured manifolds. Then, with $\intg_2$ coefficients, we have
$${SHM}(M_1\sharp M_2,\ga_1\cup\ga_2)\cong {SHM}(M_1,\ga_1)\otimes {SHM}(M_2,\ga_2)\otimes (\intg_2)^{2}.$$
\ethm

\brem
The same connected sum formula holds for sutured Instanton Floer homology, with $(\intg_2)^{2}$ replaced by $\mathbb{C}^{2}$. As a corollary, we also offer a new proof of the connected sum formula for framed instanton Floer homology:
$$I^{\sharp}(Y_1\sharp Y_2)\cong I^{\sharp}(Y_1)\otimes I^{\sharp}(Y_2).$$ 
This formula for framed instanton Floer homology has already been known to people. For example, see Scaduto \cite{scaduto2015instanton}.
\erem

The proof of the connected sum formula relies on computing the sutured monopole Floer homology of the balanced sutured manifold $(S^3(2),\delta^2)$, where $S^3(2)$ is obtained from $S^3$ by digging out two disjoint $3$-balls, and $\delta^2$ is the disjoint union of two simple closed curves, one on each spherical boundary component of $S^3(2)$. The computation in instanton theory was done by Baldwin and Sivek in \cite{baldwin2016instanton}, using an oriented skein relation in instanton theory, which was developed by Kronheimer and Mrowka in \cite{kronheimer2010instanton}. In this paper, we adapt those ideas to monopole theory and reach a similar result about oriented skein relations, and, as a corollary, we also obtain the sutured monopole Floer homology of $(S^3(2),\delta^2)$.

\begin{figure}[h]
\centering
\begin{overpic}[width=5.0in]{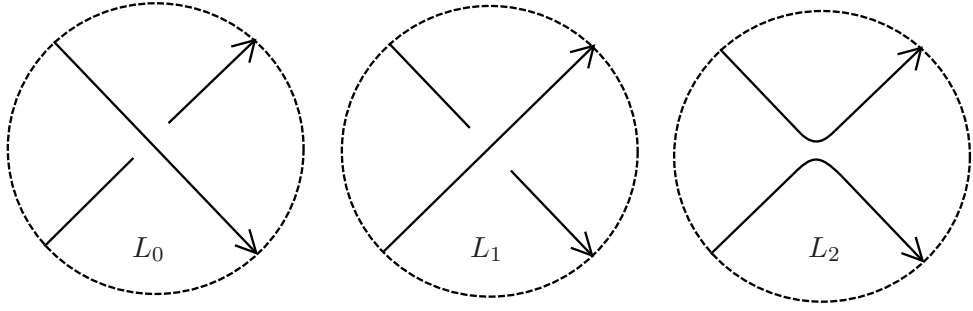}
    \put(13,5){$L_0$}
    \put(48,5){$L_1$}
    \put(83,5){$L_2$}
\end{overpic}
\vspace{0.05in}
\caption{The oriented skein relation.}\label{fig_oriented_skein_relation}
\end{figure}

\bthm\label{thm_oriented_skein_relation_introduction}
When using $\intg_2$ coefficients, there is an exact triangle for the monopole knot Floer homologies of three knots $L_0$, $L_1$ and $L_2$, which are related by an oriented skein relation as in Figure \ref{fig_oriented_skein_relation}.
\ethm

In the proof of Theorem \ref{thm_oriented_skein_relation_introduction}, another special sutured manifold $(V, \ga^4)$ arises. Here, $V$ is a framed solid torus, and the suture $\ga^4$ consists of four longitudes on $\partial{V}$. To compute the sutured monopole Floer homology of $(V,\ga^4)$, we need to obtain both a lower bound and an upper bound on its rank (fortunately, they coincide). The upper bound is obtained from by-pass exact triangles, which were introduced into the context of sutured monopole Floer homology by Baldwin and Sivek in \cite{baldwin2016contact}. The proof of the existence of by-pass exact triangles relies ultimately on the surgery exact triangle in monopole theory, which was introduced by Kronheimer, Mrowka, Ozsv\'ath, and Szab\'o in \cite{kronheimer2007monopolesandlens}. However, in \cite{kronheimer2007monopolesandlens}, they avoided the orientation issues by only working in characteristics $2$. Hence the by-pass exact triangle is only established in characteristics $2$ in Baldwin and Sivek \cite{baldwin2016contact}, and so iss our application to the computation. It is also worth mentioning that the original construction of sutured monopole Floer homology by Kronheimer and Mrowka \cite{kronheimer2010knots} does not directly apply to the case of using $\intg_2$ coefficients. It was later verified by Sivek \cite{sivek2012monopole}.

To obtain a lower bound, we construct a grading on sutured monopole Floer homology based on a properly embedded surface inside the balanced sutured manifold. The idea of such a construction originates from the proof of the decomposition theorem of suture monopole Floer homology by Kronheimer and Mrowka in \cite{kronheimer2010knots}, and was carried out in detail by Baldwin and Sivek in \cite{baldwin2018khovanov}. However, the construction by Baldwin and Sivek was restricted to the case where the properly embedded surface has a connected boundary, and it intersects the suture transversely at two points. The argument in the current paper is a naive generalization of their work, and a more systematical treatment is in Li \cite{li2019direct}.

In this paper, we prove the following theorem. 

\bthm
Let $(V,\lambda^{2n})$ be a solid torus with $2n$ longitudes as the suture. We work with $\mathbb{Q}$ coefficients, and suppose $n=2k+1$ is odd. Then, there is a grading induced by a meridian disk of $V$, and, under this grading, the sutured monopole Floer homology of $(V,\lambda^{2n})$ can be described as follows:
\begin{equation*}
{SHM}(V,\ga^{2n},i)\cong
\left\{
\begin{array}{ll}
    H_{i+k}(T^{n-1}),&-k\leq i\leq k.\\
    0&i>k~{\rm or}~i<-k.\\
\end{array}    
\right.
\end{equation*}

The same conclusion also holds for sutured instanton Floer homology with $\mathbb{C}$ coefficients.
\ethm

\brem
As commented by Yi Xie, for sutured instanton Floer homology and for odd $n$, the representation variety of a suitable closure of $(V,\ga^{2n})$ is precisely the $(n-1)$-dimensional torus $T^{n-1}$. 
\erem

As we will explain more in Subsection \ref{subsec_computing}, the following question might be interesting:

\begin{quest}
    Is ${SHM}(V,\ga^{2n})$ fully generated by the contact elements of some tight contact structures on $(V,\ga^{2n})$?
\end{quest}

{\bf Acknowledgements.} This material is based upon work supported by the National Science Foundation under Grant No. 1808794. I would like to thank my advisor Tomasz Mrowka for his enormous help, and I would like to thank Jianfeng Lin, Langte Ma, and Yi Xie for inspiring conversations. I would also like to thank the referee for many helpful comments.

\section{Preliminaries}
\subsection{Sutured monopole Floer homology}\label{subsec_sutured_monopoles}
The definitions and notations herein are consistent with the author's previous paper \cite{li2018gluing}. For more details, readers are referred to that paper. We start with the definition of balanced sutured manifolds.

\bdefn
Suppose $M$ is a compact oriented $3$-manifold with a non-empty boundary, and $\ga$ is a collection of disjoint oriented simple closed curves on $\partial{M}$ so that the following is true.

(1) $M$ has no closed components, and any component of $\partial{M}$ contains at least one component of $\ga$.

(2) We require that the orientation on $\ga$ induces an orientation on the surface $\partial{M}\backslash\ga$. This orientation, once it exists, is unique and is called the {\it canonical orientation}.

(3) Let $A(\ga)=\ga\times[-1,1]\subset \partial M$ be an annular neighborhood of $\ga\subset\partial{M}$ and let $R(\ga)=\partial{M}\backslash{\rm int}(A(\ga))$. Following, let $R_{+}(\ga)$ be the part of $R(\ga)$ where the canonical orientation coincides with the boundary orientation induced by $M$, and let $R_{-}(\ga)=R(\ga)\backslash R_{+}(\ga)$. Then, we require that
$$\chi(R_+(\ga))=\chi(R_-(\ga)).$$
The pair $(M,\ga)$ is called a {\it balanced sutured manifold}.
\edefn

To define sutured monopole Floer homology, we need to construct a closed $3$-manifold out of the sutured data. Suppose that $(M,\ga)$ is a balanced sutured manifold. Let $T$ be a compact connected oriented surface so that the following is true. 

(1) There exists an orientation reversing diffeomorphism $f:\partial{T}\ra\ga$.

(2) $T$ contains a simple closed curve $c$ that represents a non-trivial class in $H_1(T)$. 

Let
$$\widetilde{M}=M\mathop{\cup}_{f\times id} T\times[-1,1],$$
and suppose the two oriented boundary components of $\widetilde{M}$ are
$$\partial{\widetilde{M}}=R_{+}\cup R_{-}.$$
We know that $c\times\{\pm1\}\subset R_{\pm}$ is non-separating by assumption. Let 
$$h:R_+\ra R_-$$
be an orientation preserving diffeomorphism so that
$$h(c\times\{1\})=c\times\{-1\}.$$
We can use $h$ to glue the two boundary components of $\widetilde{M}$ together to obtain a closed $3$-manifold $Y$. Alternatively, we can define
$$Y=\widetilde{M}\mathop{\cup}_{id\times\{-1\}\cup h\times\{1\}} R_+\times[-1,1].$$
Let $R=R_+\times\{0\}\subset Y$.

\bdefn\label{defn_closure}
The pair $(Y,R)$ is called a {\it closure} of $(M,\ga)$. The choices $T,f,c$, and $h$ are called the {\it auxiliary data}. In particular, $T$ is called an {\it auxiliary surface}. Pick $\eta$ to be a non-separating simple closed curve on $R$, and let
$$\mathfrak{S}(Y|R)=\{\mathfrak{s}~{\rm spin}^c~{\rm structures}~Y|~c_1(\mathfrak{s})[R]=2g(R)-2\}.$$
Then, define the sutured monopole Floer homology of $(M,\ga)$ to be
$$SHM(M,\ga)=HM(Y|R;\Gamma_{\eta})=\bigoplus_{\mathfrak{s}\in\mathfrak{S}(Y|R)}\widecheck{HM}_{\bullet}(Y,\mathfrak{s};\Gamma_{\eta}).$$
\edefn

\brem\label{rem_choice_of_coefficient_ring}
The choice of the simple closed curve $\eta$ has more restrictions in Kronheimer and Mrowka's original construction in \cite{kronheimer2010knots}. The fact that we could pick any non-separating simple close curve is a simple application of the Floer excision introduced by Kronheimer and Mrowka in \cite{kronheimer2010knots}. For two different choices of the curves $\eta$, the corresponding sutured monopole Floer homologies are isomorphic to each other (and actually canonically isomorphic, as shown by Baldwin and Sivek \cite{baldwin2015naturality}). 

The curve $\eta$ may be absent, when it is convenient to use $\intg$ or $\intg_2$ coefficients. In general, when $\eta$ does exist, we will use the Novikov ring $\mathcal{R}$ or other suitable rings to construct local coefficients. Also, when defining the contact elements, the choice of the simple closed curve $\eta\subset R$ is more restrictive. For more details, readers are referred to Kronheimer and Mrowka \cite{kronheimer2010knots}, Baldwin and Sivek \cite{baldwin2016contact} and Sivek \cite{sivek2012monopole}.
\erem

The fact that sutured monopole Floer homology is well defined is proved by Kronheimer and Mrowka in \cite{kronheimer2010knots}.
\bthm
The isomorphism class of ${SHM}(M,\ga)$ is independent of all the auxiliary data and the choice of the curve $\eta$. In that way, it serves as an invariant of the balanced sutured manifold $(M,\ga)$.
\ethm

Floer excisions are used repeatedly in the paper, and, therefore, we would like to present them here for further references. Floer excisions in the context of monopole Floer theory were first introduced by Kronheimer and Mrowka \cite{kronheimer2010knots}.

Suppose $Y_1$, $Y_2$ are two closed connected oriented $3$-manifolds, and suppose that, for $i=1,2$, there is an oriented closed surface $R_i\subset Y_i$ and an oriented torus $T_i\subset Y_i$ so that $R_i\cap T_i=c_i$. Here, $c_i$ is a simple closed curve so that there is another simple closed curve $\eta_i\subset R_i$ intersecting $c_i$ transversely once. We can cut $Y_i$ along $T_i$ to get a manifold-with-boundary $\widetilde{Y}_i$ so that
$$\partial\widetilde{Y}_i=T_{i,+}\cup T_{i,-}.$$
Here, $T_{i,\pm}$ are parallel copies of $T_i$.  Let $c_{i,\pm}\subset T_{i,\pm}$ be parallel copies of $c_i$. Note that the hypothesis that $\eta_i$ intersects $c_i$ transversely once implies that $\eta$ also intersects $T_i$ transversely once. As a result, $T_i$ represents a non-trivial homology class in $H_2(Y_i)$, and $\widetilde{Y}$ is connected. Pick an orientation preserving diffeomorphism
$$h:T_{1,+}\ra T_{2,-},$$
so that
$$h(c_{1_+})=c_{2,-},~{\rm and~}h(\eta_1\cap c_{1,+})=\eta_2\cap c_{2,-}.$$
Then, we can use $h$ to glue $\widetilde{Y}_{1}$ and $\widetilde{Y}_2$ together to get an oriented connected $3$-manifold $Y$ together with an oriented connected surface $R$, which is obtained by cutting and re-gluing $R_1$ and $R_2$ together. Also, $\eta_1$ and $\eta_2$ are cut and re-glued together to form a new simple closed curve $\eta\subset R$.

\begin{figure}[h]
\centering
\begin{overpic}[width=5.0in]{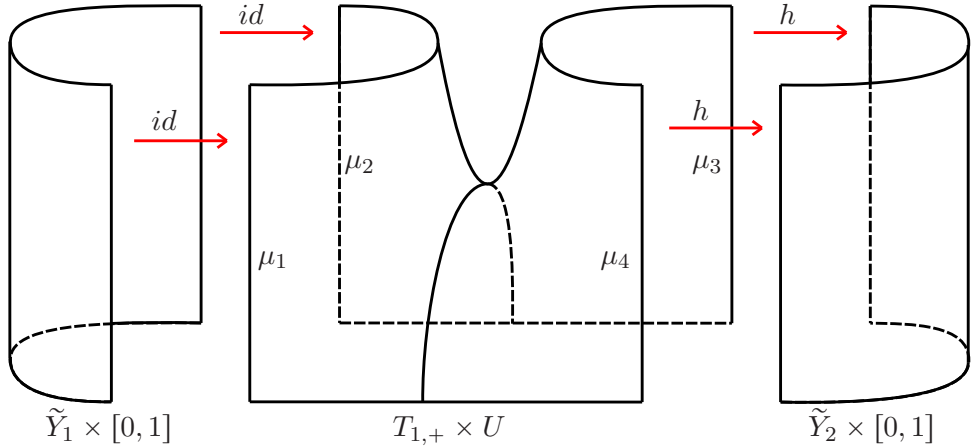}
    \put(4,-3){$\widetilde{Y}_1\times[0,1]$}
    \put(40,-3){$T_{1,+}\times U$}
    \put(83,-3){$\widetilde{Y}_2\times[0,1]$}
    \put(15,29){$id$}
    \put(24,40){$id$}
    \put(71,30){$h$}
    \put(80,40){$h$}
    \put(35,25){$\mu_2$}
    \put(26,15){$\mu_1$}
    \put(71,25){$\mu_3$}
    \put(61.5,15){$\mu_4$}
\end{overpic}
\vspace{0.05in}
\caption{Gluing three parts together to get $W$. The middle part is $T_{1,+}\times U$, while the $T_{1,+}$ directions shrink to a point in the figure.}\label{fig_excision_cobordism}
\end{figure}

Next, we construct a cobordism from $Y_1\sqcup Y_2$ to $Y$ as follows: Let $U$ be the surface as depicted in the middle part of Figure \ref{fig_excision_cobordism}, and let $\mu_1,\mu_2,\mu_3,$ and $\mu_4$ be the four vertical arcs that are part of the boundary of $U$. Suppose that each $\mu_i$ is identified with the interval $[0,1]$. 

Then, let
$$W=(\widetilde{Y}_1\times[0,1])\mathop{\cup}_{\phi}(T_{1,+}\times U)\mathop{\cup}_{\psi}(\widetilde{Y}_2\times[0,1])$$
be the $4$-manifold obtained by gluing three pieces together. Here,
$$\phi=(id\cup id)\times id:(T_{1,+}\cup T_{1,-})\times[0,1]\ra T_{1,+}\times(u_1\cup u_2)$$
and
$$\psi=(h\cup h)\times id:T_{1,+}\times(u_3\cup u_4)\ra (T_{2,+}\cup T_{2,-})\times[0,1]$$
are the gluing maps. Let $R_W=R_1\cup R_2\cup R$, and let
$$\nu=((\eta_1\cap \widetilde{Y}_1)\times[0,1])\mathop{\cup}_{\phi}((\eta_1\cap c_{1,+})\times U)\mathop{\cup}_{\psi}((\eta_2\cap \widetilde{Y}_2)\times[0,1]).$$
See Figure \ref{fig_excision_cobordism}. Then, we can define a map
\begin{equation}\label{eq_Floer_excision_map}
F=\widecheck{HM}(W|R_W;\Gamma_{\nu}):{HM}(Y_1\sqcup Y_2|R_1\cup R_2;\Gamma_{\eta_1\cup\eta_2})\ra {HM}(Y|R;\Gamma_{\eta}).
\end{equation}
In \cite{kronheimer2010knots}, Kronheimer and Mrowka proved the following theorem.
\bthm\label{thm_floer_excision}
The map $F$ in (\ref{eq_Floer_excision_map}) is an isomorphism.
\ethm

\brem
In the rest of the paper, when the choices of the surface and the local coefficients are clear in context, we will omit them from the notation and simply write:
$$\widecheck{HM}(W):{HM}(Y_1\sqcup Y_2|R_1\cup R_2)\ra{HM}(Y|R).$$
\erem

\subsection{Arc configurations and contact elements}\label{subsec_contact_element}
In this subsection, we review the construction of contact elements in sutured monopole Floer homology, in Baldwin and Sivek \cite{baldwin2016contact}.

\bdefn\label{defn_compatible_contact_structure}
Suppose $(M,\ga)$ is a balanced sutured manifold. A contact structure $\xi$ on $M$ is called {\it compatible} if $\partial{M}$ is convex and $\ga$ is (isotopic to) the dividing set on $\partial{M}$.
\edefn

\bdefn\label{defn_arc_configuration}
Suppose $T$ is a connected compact oriented surface-with-boundary. An {\it arc configuration} $\mathcal{A}$ on $T$ consists of the following data.

(1) A finite collection of pairwise disjoint simple closed curves $\{c_1,...c_m\}$ so that for any $j$, $[c_j]\neq0\in H_1(T).$

(2) A finite collection of pairwise disjoint simple arcs $\{a_1,...,a_n\}$ so that the following is true.

\quad (a) For any $i,j$, ${\rm int}(a_i)\cap c_j=\emptyset$, and ${\rm int}(a_i)\cap \partial{T}=\emptyset$.

\quad (b) For any $i$, one end point of $a_i$ lies on $\partial{T}$ and the other on some $c_j$.

\quad (c) Each boundary component of $T$ has a non-trivial intersection with some $a_i$.

\begin{figure}[h]
\centering
\begin{overpic}[width=4.0in]{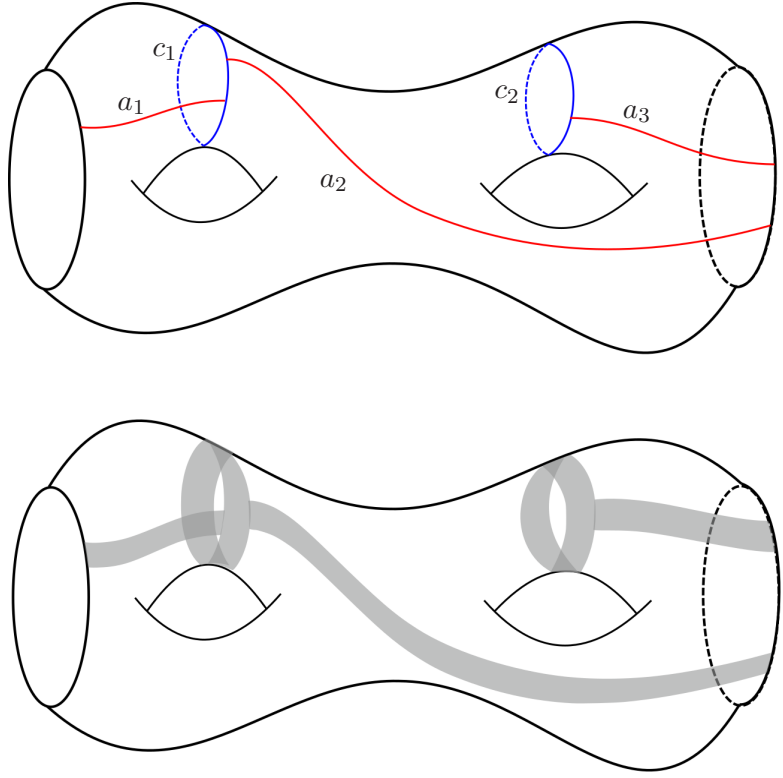}
    \put(14,86){$a_1$}
    \put(18.5,93){$c_1$}
    \put(40,76){$a_2$}
    \put(79,85){$a_3$}
    \put(62.5,87.5){$c_2$}
\end{overpic}
\vspace{0.05in}
\caption{Top, an arc configuration on $T$. Bottom, the shaded region corresponds to the negative region on $T\times\{t\}\subset T\times[-1,1]$, with respect to the contact structure induced by the arc configuration. Its boundary is the dividing set on $T\times\{t\}$.}\label{fig_arc_configuration}
\end{figure}

See Figure \ref{fig_arc_configuration} for an example of the arc configuration $\mathcal{A}$. An arc configuration is called {\it reduced} if there is only one simple closed curve in (1). 
\edefn

Let $(M,\ga)$ be a balanced sutured manifold equipped with a compatible contact structure $\xi$. Suppose $T$ is a connected auxiliary surface for $(M,\ga)$ and $\mathcal{A}$ is a reduced arc configuration on $T$. One can constructe a suitable contact structure $\widetilde{\xi}$ on
$$\widetilde{M}=M\cup T\times[-1,1]$$
as follows: First, the arc configuration $\mathcal{A}$ gives rise to an $[-1,1]$-invariant contact structure on $T\times[-1,1]$. The negative region on any piece $T\times\{t\}$ is shown in Figure \ref{fig_arc_configuration}. Then, perturbe the contact structure on $M$ in a neighborhood of $\ga\subset M$ so that the dividing set on $A(\ga)$ can be identified with the one on $\partial{T}\times[-1,1]$. Thus, there is an orientation reversing diffeomorphism $f:\partial{T}\times[-1,1]\ra A(\ga)$, which also identifies the contact structures. We can use $f$ to glue $T\times [-1,1]$ to $M$, and, after rounding the corners, the desired contact structure $\widetilde{\xi}$ on $\widetilde{M}$ is constructed. Suppose that
$$\partial{\widetilde{M}}=R_+\cup R_-,$$
then $R_{\pm}$ are convex and the dividing set on $R_+$ or $R_-$ consists of two parallel non-separating simple closed curves. Finally, choose a gluing diffeomorphism $h:R_+\ra R_-$, which preserves the contact structures, to construct a closure $(Y,R)$, equipped with a contact structure $\bar{\xi}$. Under $\xi$, $R$ is convex and the negative region on $R$ is just an annulus. Also, choose a simple closed curve $\eta\subset R$, which intersects each component of the dividing set transversely once, to support the local coefficients. From the construction, we know that
$$c_1(\bar{\xi})[R]=2-2g(R),$$
and, thus, by Kronheimer, Mrowka, Ozsv\'ath, and Szab\'o \cite{kronheimer2007monopolesandlens}, there is a contact element
$$\phi_{\bar{\xi}}\in \widecheck{HM}(-Y,\mathfrak{s}_{\bar{\xi}};\Gamma_{-\eta})\subset{SHM}(-M,-\ga).$$

\bthm[Baldwin and Sivek \cite{baldwin2016contact}]
The contact element
$$\phi_{\bar{\xi}}\in{SHM}(-M,-\ga),$$
which is constructed by only using reduced arc configurations, is independent of the choices made in the construction and, thus, serves as a well-defined invariant of $\xi$ on $(M,\ga)$.
\ethm

\brem
In \cite{baldwin2016contact}, Baldwin and Sivek only used reduced arc configurations to construct contact elements. However, part of the construction can be made with a general arc configuration, as defined in Definition \ref{defn_arc_configuration}. When using an arc configuration that is not necessarily reduced, the dividing set on $R_{\pm}$ consists of $m$ many pairs of parallel non-separating simple closed curves. Here, $m$ is the number of simple closed curves in the arc configuration. However, in this case, the diffeomorphism $h$, which preserves the contact structures on $R_{\pm}$, may not always exist (as it must identify the dividing sets). The reason we want to introducing non-reduced arc configurations is that, as we see in later sections, a general arc configuration will arise when performing Floer excisions, and in that special case, the diffeomorphism $h$ can indeed be chosen.
\erem

To conclude the current section, we introduce the definition of contact handle attachments for further references in Subsection \ref{subsec_connected_sum}.

\bdefn
Suppose $(M,\ga)$ is a balanced sutured manifold equipped with a compatible contact structure. A {\it contact handle} attached to $(M,\ga)$ is a quadruple $h=(\phi,S,D^3,\delta)$ so that the following is true.

(1) $D^3$ is a $3$-ball equipped with the standard tight contact structure and $\delta$ is the dividing set on $\partial{D}^3$.

(2) $S\subset \partial{D}^3$ is a compact submanifold and $\phi:S\ra \partial{M}$ is an embedding so that $\phi(S\cap \delta)\subset \ga$. $S$ has different descriptions due to the index of the contact handle.

\quad (a) In the index $0$ case, $S=\emptyset$.

\quad (b) In the index $1$ case, $S$ is a disjoint union of two disks, and each disk intersects $\delta$ in an arc.

\quad (c) In the index $2$ case, $S$ is an annulus intersecting $\delta$ in two arcs. Also, we require that each component of $\partial{S}$ intersects each arc transversely once.

\quad (d) In the index $3$ case, $S=\partial{D}^3$.
\edefn

\section{Contact element and excision}
Suppose, for $i=1,2$, $(M_i,\ga_i)$ is a balanced sutured manifold. Suppose further that $(T_i,f_i,c_i,h_i)$ is the auxiliary data for constructing a closure $(Y_i,R_i)$ of $(M_i,\ga_i)$, as in Definition \ref{defn_closure}. According to the construction of closures in Section \ref{subsec_sutured_monopoles}, $R_i$ contains a curve corresponding to the curve $c_i\subset T_i$. Abusing the notation, we also denote this curve on $R_i$ by $c_i$. We can choose a simple closed curve $\eta_i$ on $R_i$ with exactly one transverse intersection with $c_i$. 

\begin{figure}[h]
\centering
\begin{overpic}[width=4.0in]{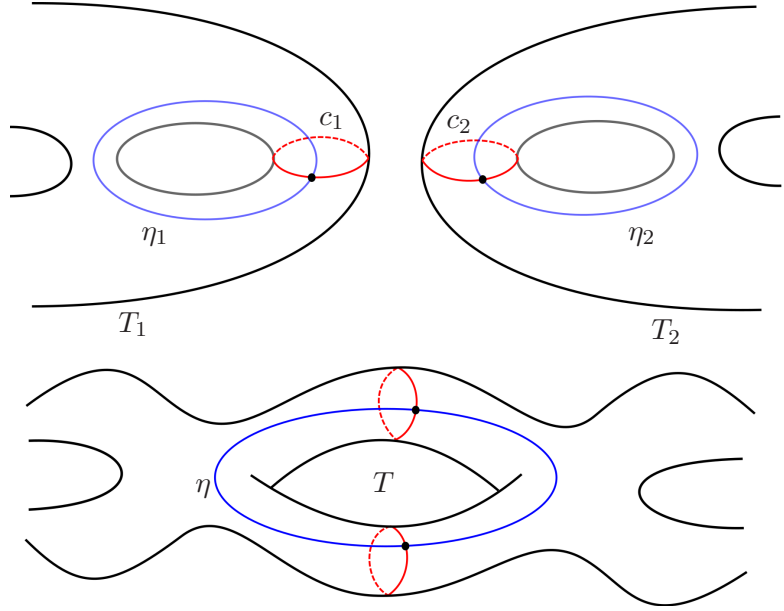}
    \put(14,36){$T_1$}
    \put(83,35){$T_2$}
    \put(47,15){$T$}
    \put(17,48){$\eta_1$}
    \put(80,48){$\eta_2$}
    \put(24,15){$\eta$}
    \put(40,63){$c_1$}
    \put(56.5,62.5){$c_2$}
\end{overpic}
\vspace{0.05in}
\caption{Top, the two auxiliary surfaces $T_1$ and $T_2$. Bottom, the connected auxiliary surface $T$.}\label{fig_cut_and_paste}
\end{figure}

Let $M=M_1\sqcup M_2$ and $\ga=\ga_1\cup\ga_2$. Then, $(M,\ga)$ is also a balanced sutured manifold. We can cut $T_i$ along $c_i$ and re-glue the newly created boundary with respect to the orientation. Then, $T_1$ and $T_2$ become a connected surface $T$ so that
$$g(T)=g(T_1)+g(T_2)-1,~\partial{T}=\partial{T_1}\cup\partial{T}_2.$$
Choose $f=f_1\cup f_2$ and $h=h_1\cup h_2$. When cutting and re-gluing along $c_1$ and $c_2$, the two curves $\eta_1$ and $\eta_2$ can also be cut and glued together to become a curve $\eta$. We can use the auxiliary data $(T,f,h,\eta)$ to close up $(M,\ga)$ and then obtain a closure $(Y,R)$ of $(M,\ga)$. See Figure \ref{fig_cut_and_paste}. As in Subsection \ref{subsec_sutured_monopoles}, we can construct a Floer excision map:
\begin{equation}\label{eq_Floer_excision_maps_2}
F:\widecheck{HM}(-(Y_1\sqcup Y_2)|-(R_1\cup R_2))\ra\widecheck{HM}(-Y|-R).
\end{equation}
We have the following theorem.

\bthm\label{thm_excision_preserves_contact_element}
Suppose that the genus of $T_1$ and $T_2$ are large enough, and, for $i=1,2$, $(M_i,\ga_i)$ is equipped with a compatible contact structure $\xi_i$. Suppose further that $t$ is obtained from $T_1$ and $T_2$ but cutting and re-gluing as described above. Then, we can find suitable reduced arc configurations $\mathcal{A}_1,\mathcal{A}_2$, and $\mathcal{A}$ on $T_1,T_2$, and $T$, respectively to construct contact structures $\bar{\xi}_1,\bar{\xi}_2$, and $\bar{\xi}$ on suitable closures $Y_1,Y_2$, and $Y$. Then, the map $F$ in (\ref{eq_Floer_excision_maps_2}) preserves contact elements:
$$F(\phi_{\bar{\xi}})\doteq \phi_{\bar{\xi}_1\cup\bar{\xi}_2}.$$
Here, $\doteq$ means equal up to multiplication by a unit.
\ethm

\begin{figure}[h]
\centering
\begin{overpic}[width=4.0in]{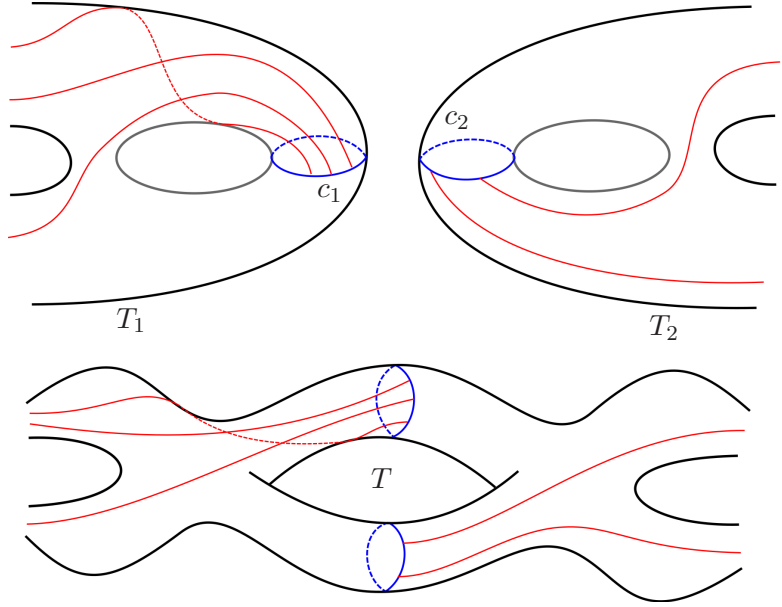}
    \put(14,36){$T_1$}
    \put(83,35){$T_2$}
    \put(47,15){$T$}
    \put(40,53){$c_1$}
    \put(56.5,62.5){$c_2$}
\end{overpic}
\vspace{0.05in}
\caption{Top, The two reduced arc configurations on $T_1$ and $T_2$. Bottom, the resulting arc configuration on $T$ from slicing. It has two simple closed curves instead of one.}\label{fig_slicing}
\end{figure}

\bpf
We first choose the special arc configurations $\mathcal{A}_1$ and $\mathcal{A}_2$. For $i=1,2$, assume that we have a reduced arc configuration $\mathcal{A}_i$ on $T_i$ so that the simple closed curve is $c_i$, and all arcs are attached to only one side of $c_i\subset T_i$. See Figure \ref{fig_slicing}.

To show that such special arc configurations $\mathcal{A}_1$ and $\mathcal{A}_2$ do exist, we work only with $(T_1,c_1)$, and the argument for $(T_2,c_2)$ is similar. Cut $T_1$ open along $c_1$ and let $c_{1,+}$ $c_{1,-}$ be the two newly created boundary components. Now, it is enough to find a set of pair-wise disjoint properly embedded arcs ${a_1,...,a_n}$ on $T_1\backslash c_1$ so that the following is true.

(1) For any $i\in\{1,...,n\}$, one end point of $a_i$ is on $\partial{T}_i$, and the other end point of $a_i$ is on $c_{1,+}$.

(2) Each component of $\partial{T}_i$ intersects with some $a_i$.

We can pick the set of arcs $\{a_1,...,a_n\}$ one by one. First, pick any $a_1$ that satisfies (1) and is non-separating on $T_1\backslash c_1$. Then, we can pick an arc $a_2$ that connects a different component of $\partial{T}_1$ to $c_{1,+}$ and is non-separating in $T_1\backslash(c_1\cup a_1)$. Keep running this process until all boundary components of $\partial{T}_1$ have been connected to $c_{1,+}$ by an arc. Note that we can make the genus of $T_1$ as large as we want, which means that it is always possible to find the desired set of arcs. To obtain the arc configuration $\mathcal{A}_1$, we glue $c_{1,+}$ to $c_{1,-}$ to recover $T_1$. Since all arcs are chosen to be attached to $c_{1,+}$, they appear on the same side of $c_1$ on $T_1$.

 Since the arcs $a_j$ are attached to the same side of $c_i$, when looking at the induced contact structure on $T_i\times[-1,1]$, the boundary of the negative region on $T_i\times\{t\}$ consists of a few arcs, whose end points are on $\partial{T_i}\times\{t\}$, and one simple closed curve, which is a parallel copy of $c_i\subset T_i$. Abusing the notation, we still use $c_i$ to denote this closed component of the boundary of the negative region on $T_i\times\{t\}$. We then pick a gluing diffeomorphism $h_i$ that identifies the contact structures on the boundary of $\widetilde{M}_i=M_i\cup T_i\times[-1,1]$ and which also preserves $c_i$.

When we extend $\xi_i$ to $\bar{\xi}_i$, which is defined on all of $Y_i$, the new contact structure $\bar{\xi}_i$ will be $S^1$-invariant in a neighborhood of $c_i$. To describe this contact structure in coordinates, let $A_i\subset T_i$ be a neighborhood of $c_i\subset T_i$. In $Y_i$, $A_i\times S^1$ is a neighborhood of $c_i\subset Y_i$. In this neighborhood, we can write the contact form as
$$\al_i=\be_i+u_i\cdot d\varphi_i,$$
where $\be_i$ is a $1$-form on $A_i$, $u_i$ is a function on $A_i$ with
$$c_i=\{p\in A_i|u_i(p)=0\},$$
and $\varphi_i$ is the coordinate for $S^1$ direction. See Geiges \cite{geiges2008introduction}. The non-degeneracy condition reads
\beq
0\neq \al_i\wedge d\al_i
=(u_i\cdot d\be_i+\be_i\wedge du_i)\wedge d\varphi_i.
\eeq
Along $c_i$, we have $\be_i\wedge du_i\neq0.$ Hence, along $c_i$, the $d\theta_i$ component of $\be$ is always non-zero. Here, $\theta_i$ is a coordinate for $c_i$, and $(u,\theta_i)$ can serve as local coordinates in a small neighborhood $[-\varepsilon,\varepsilon]\times c_i\subset A_i$. Then, the slicing operation defined by Niederkr\"uger and Wendl in \cite{niederkruger2011weak} can be described as follows: Let $L_i=c_i\times S^1$ be a pre-Lagrangian torus, (for the definition of pre-Lagrangian tori, see Ma \cite{ma2000pre}) and let $N_i=[-\varepsilon,\varepsilon]\times c_i\times S^1$ be a neighborhood of $L_i$ with the coordinates $(u_i,\theta_i,\varphi_i)$. Note the coordinate $u_i$ corresponds to $r$ in \cite{niederkruger2011weak}, and the other two coordinates are the same as in that paper. We can cut $N_i$ open along $L_i$ so that $N_i$ is cut into two parts $N_{i,+}$ and $N_{i,-}$, which correspond to where $u_i\geq0$ and $u_i\geq 0$, respectively. Then, re-glue $N_{1,+}$ to $N_{2,-}$ and $N_{1,-}$ to $N_{2,+}$ by identifying $L_1$ with $L_2$ so that $(\theta_1,\varphi_1)$ is identified with $(\theta_2,\varphi_2)$. Suppose that the resulting $3$-manifold is $Y$, then $Y$ has a distinguishing surface $R$ obtained by cutting and re-gluing $R_1$ and $R_2$ along $c_1$ and $c_2$. Recall that there is a simple closed curve $\eta_i\subset R_i$, which intersects $c_i$ transversely once. After a suitable isotopy, we can assume that, under the above identification of $L_1$ with $L_2$, we can also identify $\eta_1\cap c_1$ with $\eta_2\cap c_2$. Hence, $\eta_1$ and $\eta_2$ are also cut and re-glued to become a curve $\eta\subset R$. This is exactly the same procedure of performing a Floer excision along the tori $L_1$ and $L_2$. Thus, $(Y,R)$ is a closure of $(M_1\sqcup M_2,\ga_1\cup\ga_2)$. Also, by Theorem \ref{thm_floer_excision}, there is an isomorphism
\begin{equation*}
    F: HM(-(Y_1\sqcup Y_2)|-(R_1\cup R_2);\Gamma_{-(\eta_1\cup\eta_2)})\ra HM(-Y|-R;\Gamma_{-\eta}).
\end{equation*}

\begin{figure}[h]
\centering
\begin{overpic}[width=4.0in]{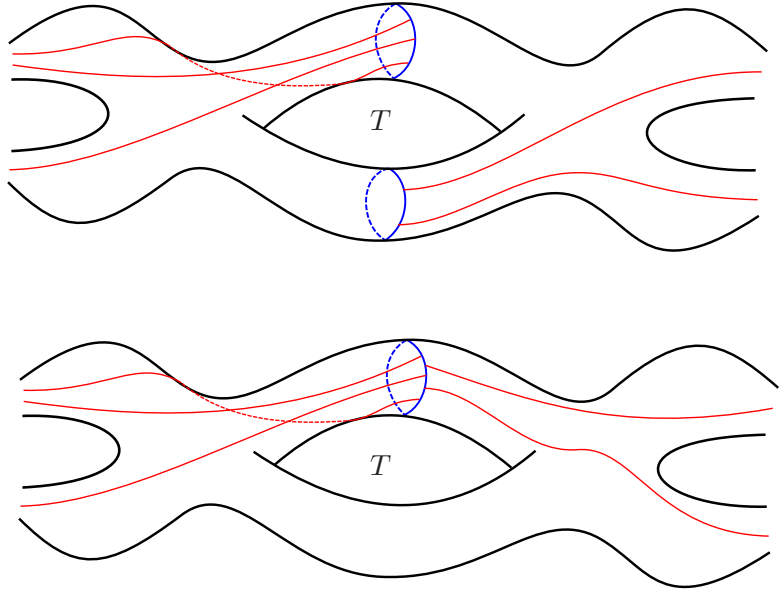}
    \put(47,60){$T$}
    \put(47,15){$T$}
\end{overpic}
\vspace{0.05in}
\caption{Top, the arc configuration on $T$ obtained from slicing. Bottom, the reduced arc configuration after merging the two simple closed curves.}\label{fig_slicing2}
\end{figure}

The process of slicing also cuts and re-glues the contact structures $\bar{\xi}_i$ on $Y_i$ to obtain a contact structure $\bar{\xi}'$ on $Y$, as explained in \cite{niederkruger2011weak}. The contact structure $\bar{\xi}'$, however, arises from an arc configuration $\mathcal{A}'$ that is not reduced in the sense of Definition \ref{defn_arc_configuration}. This is because, with respect to $\bar{\xi}'$, the dividing set on $R$ consists of two pairs of parallel non-separating simple closed curves rather than just one pair. See Figure \ref{fig_slicing}. Let $\bar{\xi}$ be a contact structure on $Y$, which is obtained by extending $\xi_i$ on $M_i$ using a reduced arc configuration $\mathcal{A}$. Here, $\mathcal{A}$ is obtained by “merging” the two simple closed curves of $\mathcal{A}'$ into one as depicted in Figure \ref{fig_slicing2}. The proof of Theorem \ref{thm_excision_preserves_contact_element} is clearly the combination of the following two lemmas.

\epf

\blem\label{lem_identify_contact_elements_for_different_contact_structure}
If the genus of $T_1$ and $T_2$ are large enough, then the spin${}^c$ structures associated to $\bar{\xi}$ and $\bar{\xi}'$ are the same. Furthermore, if we denote that spin${}^c$ structure by $\mathfrak{s}_0$, then we have
$$\phi_{\bar{\xi}}\doteq\phi_{\bar{\xi}'}\in \widecheck{HM}(-Y,\mathfrak{s}_0;\Gamma_{-\eta}).$$
\elem

\blem\label{lem_identify_contact_elements_through_excision}
If the genus of $T_1$ and $T_2$ are large enough, then we have
$$F(\phi_{\bar{\xi}_1}\cup\phi_{\bar{\xi}_2})\doteq \phi_{\bar{\xi}'}.$$
\elem

To prove the above two lemmas, we first need some preliminaries.

\blem[Baldwin and Sivek \cite{baldwin2016contact}]\label{lem_partial_open_book_decomposition}
Suppose $(M,\ga)$ is a balanced sutured manifold equipped with a compatible contact structure $\xi$ and $T$ is a connected auxiliary surface with a large enough genus. Suppose further that we use an arc configuration (not necessarily reduced) on $T$ to extend $\xi$ to a contact structure $\bar{\xi}$ on a suitable closure $(Y,R)$ of $(M,\ga)$. Then, there is a contact structure $\xi_R$ on $R\times S^1$ and a set of pair-wise disjoint simple closed curves $\{\al_1,...,\al_n\}$ so that the following is true.

(1) The contact structure $\xi_R$ is $S^1$-invariant so that each $R\times \{t\}$ is convex with the dividing set being some pairs of parallel non-separating simple closed curves.

(2) Each $\al_i$ is Legendrian and is disjoint from the pre-Lagrangian tori of the form
$$({\rm Dividing}~{\rm set}~{\rm on}~R)\times S^1.$$

(3) The result of performing $+1$ contact surgeries along all $a_i\subset R\times S^1$ is contactomorphic to $(Y,\bar{\xi})$.
\elem

\blem[Niederkr\"uger and Wendl \cite{niederkruger2011weak}]\label{lem_weakly_fillable}
Suppose $R$ is the surface as described in Lemma \ref{lem_partial_open_book_decomposition}, and $\xi_R$ is an $S^1$-invariant contact structure on $R\times S^1$ so that each $R\times{t}$ is convex with the dividing set being a few pairs of non-separating simple closed curves. Suppose further that there is a curve $\eta\subset R$, which intersects every component of the dividing set transversely once. Then, $(R\times S^1,\xi_R)$ is weakly fillable by some $(W,\omega)$ so that $\eta$ is dual to $\omega|_{R\times S^1}$ up to a (non-zero) scalar.
\elem

\blem[Kronheimer, Mrowka, Ozsv\'ath, and Szab\'o \cite{kronheimer2007monopolesandlens}]\label{lem_non_vanishing_contact_invariant}
In Lemma \ref{lem_weakly_fillable}, the contact element
$$\phi_{\xi_{R}}\in\widecheck{HM}(-R\times S^1,\mathfrak{s}_{\xi_R};\Gamma_{-\eta})$$
is primitive.
\elem

\blem[Kronheimer and Mrowka \cite{kronheimer2010knots}]\label{lem_unique_spin_c_structure}
Suppose $R$ is the surface as in Lemma \ref{lem_partial_open_book_decomposition}. Then, there is a unique spin${}^c$ structure $\mathfrak{s}_0$ on $R\times S^1$ so that the following is true.

(1) We have $c_1(\mathfrak{s}_0)[R]=2-2g.$

(2) The monopole Floer homology of $\widecheck{HM}(-R\times S^1,\mathfrak{s}_0;\Gamma_{-\eta})$ is non-zero.
Furthermore, we have 
$$\widecheck{HM}(-R\times S^1,\mathfrak{s}_0;\Gamma_{-\eta})\cong\mathcal{R}.$$
Here, $\mathcal{R}$ is the base ring we use to construct local coefficients, as in Remark \ref{rem_choice_of_coefficient_ring}.
\elem

\blem[Baldwin and Sivek \cite{baldwin2016contact}]\label{lem_surgery_preserve_contact_element} 
Suppose, for $i=1,2$, $Y_i$ is a closed oriented $3$-manifold equipped with contact structure $\xi_i$. Suppose further that $(Y_2,\xi_2)$ is obtained from $(Y_1,\xi_1)$ by performing a contact $+1$ surgery along a Legendrian curve. Then, there is a cobordism $W$, from $Y_1$ to $Y_2$, obtained from $Y_1\times[0,1]$ by attaching a $2$-handle with a suitable framing. Suppose $\eta_1$ is a $1$-cycle in $Y_1$ supporting local coefficients and is disjoint from the Legendrian curve along which we perform the Dehn surgery. Thus, $\eta_1$ remains in $Y_2$, and we call it $\eta_2$. Then, the map
$$\widecheck{HM}(-W):\widecheck{HM}(-Y_1,\mathfrak{s}_{\xi_1};\Gamma_{-\eta_1})\ra\widecheck{HM}(-Y_2,\mathfrak{s}_{\xi_2};\Gamma_{-\eta_2})$$
preserves the contact elements (up to multiplication by a unit).
\elem

\brem
Lemma \ref{lem_surgery_preserve_contact_element} is stated in Baldwin and Sivek \cite{baldwin2016contact}, as a corollary to a result in Hutchings and Taubes \cite{hutchings2013proof}.
\erem

\bpf[Proof of lemma \ref{lem_identify_contact_elements_for_different_contact_structure}]
As in the settings of Theorem \ref{thm_excision_preserves_contact_element}, $\bar{\xi}$ and $\bar{\xi}'$ are contact structures on $Y$, which are obtained from the contact structures $\xi_1\cup\xi_2$ on $(M_1,\ga_1)\cup (M_2,\ga_2)$ and some particular arc configurations $\mathcal{A}$ and $\mathcal{A}'$ on $T$. From Lemma \ref{lem_partial_open_book_decomposition}, we know that there are contact structures $\xi_R$ and $\xi_R'$ on $R\times S^1$ and a set of pair-wise disjoint curves $\al_1,...,\al_n\subset R\times S^1$ so that the following is true.

(1) Both $\xi_R$ and $\xi_R'$ are $S^1$-invariant, and any $R\times\{t\}$, for $t\in S^1$, is convex.

(2) We have $\xi_R=\xi_R'$ near a neighborhood of each $\al_i$.

(3) All $\al_i$ are disjoint from the pre-Lagrangian tori of the form
$$({\rm Dividing}~{\rm set}~{\rm on}~R)\times S^1,$$
for the dividing sets with respect to both $\xi_R$ and $\xi_R'$.

(4) If we perform contact $+1$ surgery along all of $\al_i$, then $(R\times S^1,\xi_R)$ [or $(R\times S^1,\xi_R')$] will become a contact manifold that is contactomorphic to $(Y,\bar{\xi})$ [or $(Y,\bar{\xi}')$].

Condition (2) relies on the proof of Lemma \ref{lem_partial_open_book_decomposition} (of the current paper) in \cite{baldwin2016contact}. The essential reason is that $\bar{\xi}$ and $\bar{\xi}'$ are only different in the part of $Y$ that comes from gluing auxiliary surfaces, while the curves $\al_i$ are contained in the interior of the original balanced sutured manifold.

Via Lemmas \ref{lem_weakly_fillable} and \ref{lem_non_vanishing_contact_invariant}, we know that the contact invariants $\phi_{\xi_R}$ and $\phi_{\xi_R'}$ are both primitive in the same monopole Floer homology. Then, Lemma \ref{lem_unique_spin_c_structure} makes sure that $\xi_R$ and $\xi_R'$ correspond to the same spin${}^c$ structure $\mathfrak{s}_0$ on $R\times S^1$ (because there is only one candidate for the spin${}^c$ structures). Thus, we have
\begin{equation}\label{eq_identify_contact_element_on_R_times_S_1}
\phi_{\xi_R}\doteq\phi_{\xi_R'}\in \widecheck{HM}(-R\times S^1,\mathfrak{s}_0;\Gamma_{-\eta}),
\end{equation}
for a suitable choice of local coefficients.

The surgery description in condition (4) makes sure that, on $Y$, $\bar{\xi}$ and $\bar{\xi}'$ correspond to the same spin${}^c$ structure. This fact, together with Lemma \ref{lem_surgery_preserve_contact_element} and equality (\ref{eq_identify_contact_element_on_R_times_S_1}), implies Lemma \ref{lem_identify_contact_elements_for_different_contact_structure}. 
\epf

\bpf[Proof of Lemma \ref{lem_identify_contact_elements_through_excision}.]
First, by applying Lemma \ref{lem_partial_open_book_decomposition} to $(Y_i,\bar{\xi}_i)$, for $i=1,2$, we obtain a contact structure $\xi_{R_i}$ on $R_i\times S^1$ and a set of Legendrian curves $\{\al_{i,1},...,\al_{i,n_i}\}$ satisfying the conclusion of the lemma. In particular, if we perform contact $+1$ surgery along all of $\al_{i,j}$, we will obtain $(Y_i,\bar{\xi}_i)$. If we pick a suitable connected component $c_i$ of the dividing set on $R_i\times{t}$ and perform the slicing operation on $R_1\times S^1$ and $R_2\times S^1$, along the two pre-Lagrangian tori $c_1\times S^1$ and $c_2\times S^1$, then the result is the $3$-manifold $R\times S^1$ with the contact structure $\xi_{R}'$, as in the proof of Lemma \ref{lem_identify_contact_elements_for_different_contact_structure}. Also, the two sets of curves $\{\al_{1,1},...,\al_{1,n_1}\}$ and $\{\al_{2,1},...,\al_{2,n_2}\}$ together form the set of curves $\{\al_1,...,\al_n\}$, as in the proof of Lemma \ref{lem_identify_contact_elements_for_different_contact_structure}. There is a cobordism associated to the slicing operation or, equivalently, performing a Floer excision on $R_1\times S^1$ and $R_2\times S^1$. We call this cobordism $W_e$, and it is from $(R_{1}\times S^1)\sqcup (R_2\times S^1)$ to $R\times S^1$. There is a second cobordism $W_s$, associated to the surgeries along all of $\al_i$, as in Lemma \ref{lem_surgery_preserve_contact_element}, from $R\times S^1$ to $Y$. Finally, there is a third one, $W_F$, corresponding to the map $F$ (which is also obtained from a Floer excision again), from $Y$ to $Y_1\sqcup Y_2$.

As usual, we choose suitable surfaces and local coefficients to make the cobordism map precise, but we omit them from the notation. The map $HM(-W_e)$ preserves contact elements because it is an isomorphism between two copies of $\mathcal{R}$, and the two contact elements are both units in the corresponding copy of $\mathcal{R}$. Furthermore, the map $HM(-W_s)$ preserves contact elements by Lemma \ref{lem_surgery_preserve_contact_element}. So, if we could prove that the composition $HM(-(W_e\cup W_s\cup W_F))$ preserves the contact elements, then so does $HM(-W_{F})=F$, and, thus, Lemma \ref{lem_identify_contact_elements_through_excision} follows. 

\begin{figure}[h]
\centering
\begin{overpic}[width=5.0in]{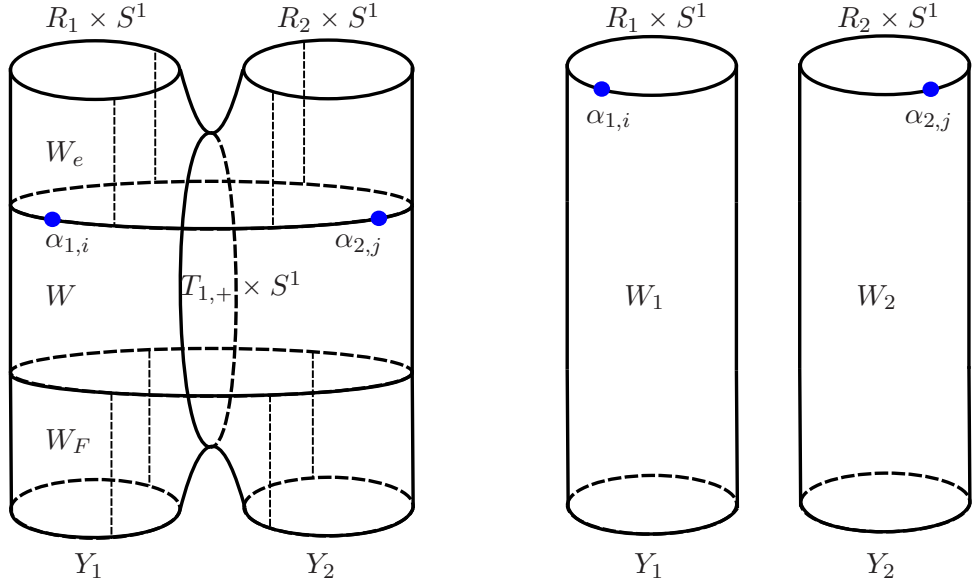}
    \put(4,54){$R_1\times S^1$}
    \put(28,54){$R_2\times S^1$}
    \put(62,54){$R_1\times S^1$}
    \put(86,54){$R_2\times S^1$}
    \put(4,31){$\al_{1,i}$}
    \put(4,25){$W$}
    \put(4,40){$W_e$}
    \put(4,10){$W_F$}
    \put(64,25){$W_1$}
    \put(88,25){$W_2$}
    \put(34,31){$\al_{2,j}$}
    \put(18,26){$T_{1,+}\times S^1$}
    \put(7,-3){$Y_1$}
    \put(31,-3){$Y_2$}
    \put(65,-3){$Y_1$}
    \put(89,-3){$Y_2$}
    \put(60,44){$\al_{1,i}$}
    \put(93,44){$\al_{2,j}$}
\end{overpic}
\vspace{0.05in}
\caption{Left, the union of the three cobordisms, cut along the $3$-torus $T_{1,+}\times S^1$. Right, the two disjoint cobordisms resulting from the cutting and pasting.}\label{fig_split_apart}
\end{figure}

To show that $HM(-(W_e\cup W_s\cup W_F))$ preserves contact elements, we observe that, when we cut the cobordism $W_e\cup W_s\cup W_F$ open along $T_{1,+}\times S^1$ and glue back two copies of $T_{1,+}\times D^2$, the result is the disjoint union of two cobordisms, which we call $W_1$ and $W_2$, respectively. See Figure \ref{fig_split_apart}. For $i=1,2$, $W_i$ is from $R_i\times S^1$ to $Y_i$ and is associated to the surgeries along $\al_{i,j}$, as in Lemma \ref{lem_surgery_preserve_contact_element}. Hence, by that lemma, $HM(-(W_1\cup W_2))$ would preserve contact elements. Finally, by Lemma 2.10 in \cite{kronheimer2010knots}, we know that
$$HM(-(W_e\cup W_s\cup W_F))\doteq HM(-(W_1\cup W_2)).$$
Thus, we conclude the proof of Lemma \ref{lem_identify_contact_elements_through_excision}.
\epf

\section{Connected sum formula}
In this section, we derive the connected sum formula for sutured monopole Floer homology. The proof of the formula relies on the computation on some special balanced sutured manifolds.

\subsection{Computing ${SHM}(V,\ga^{4k+2})$}\label{subsec_computing}

We start with the family of balanced sutured manifolds $(V,\ga^{2n})$. For $n\in\intg_+$, suppose $V=S^1\times D^2$ is a solid torus, and $\ga^{2n}\subset \partial{V}$ is the suture consisting of $2n$ longitudes (each of the form $S^1\times \{t\}$ for some $t\in\partial{D}$). To make $(M,\ga^{2n})$ a balanced sutured manifold, adjacent longitudes must be oriented oppositely.

When $n>2$, we can pick a properly embedded annulus $A$ inside $V$ so that the following is true.

(1) $\partial{A}\cap \ga^{2n}=\emptyset$.

(2) $\partial{V}\backslash\partial{A}$ has two components, and one component contains precisely three components of the suture $\ga^{2n}$ in its interior.

\begin{figure}[h]
\centering
\begin{overpic}[width=5.0in]{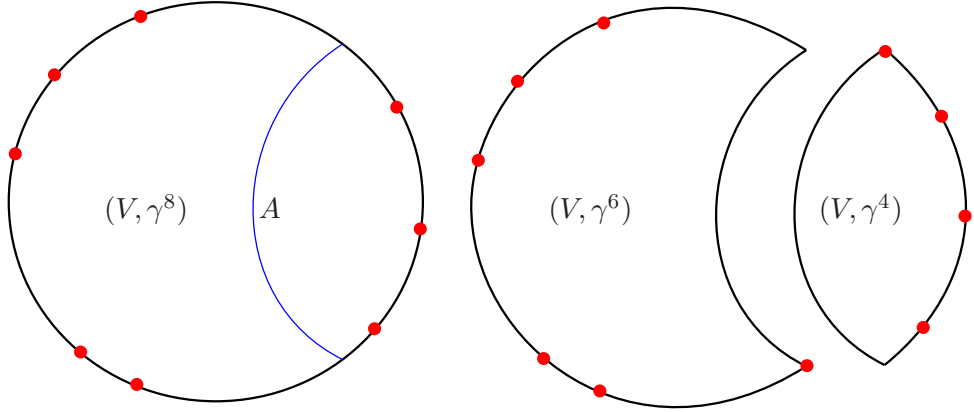}
    \put(10,20){$(V,\ga^8)$}
    \put(26,20){$A$}
    \put(56,20){$(V,\ga^6)$}
    \put(84,20){$(V,\ga^4)$}
\end{overpic}
\vspace{0.05in}
\caption{A cross section $\{t\}\times D^2$ of the solid torus$V=S^1\times D^2$. The (red) dots represent the suture and the (blue) arc inside the disk represents the annulus $A$, along which we perform the decomposition.}\label{fig_cut_along_annulus}
\end{figure}

The sutured manifold decomposition of $(V,\ga^{2n})$ along $A$ yields a balanced sutured manifold which has two components. One component is diffeomorphic to $(V,\ga^{2n-2})$, and the other is diffeomorphic to $(V,\ga^{4})$. See Figure \ref{fig_cut_along_annulus} for an example of decomposing $(V,\ga^8)$. By induction and Proposition 6.7 in Kronheimer and Mrowka \cite{kronheimer2010knots}, we know that, with $\mathbb{Q}$ coefficients and with $n\geq 2$, we have
\begin{equation}\label{eq_inductive_formula}
{SHM}(V,\ga^{2n};\mathbb{Q})\cong{SHM}(V,\ga^{4};\mathbb{Q})^{\otimes(n-1)}
\end{equation}

\blem\label{lem_four_longitudes_Z_coefficients}
When using $\intg$ coefficients, we have
$${SHM}(V,\ga^4;\intg)\cong \intg^2\oplus G_{tor},$$
where $G_{tor}$ is a (finite) torsion group without any even-torsion.
\elem

\bpf 
We first prove that the rank of the homology is precisely $2$. To get a lower bound, we use $\mathbb{Q}$ coefficients and deal with $(V,\ga^6)$. Recall that $V=S^1\times D^2$ is a solid torus. Let $t_0\in S^1$ be a point and $D=\{t_0\}\times D^2\subset V$ be a meridian disk of $V$. Note that $\partial{D}$ intersects $\ga^6$ at six points:
$$\partial{D}\cap \ga^6=\{t_0\}\times\{p_1,...,p_6\}\subset S^1\times\partial{D^2}.$$

Let $p_i$ be indexed according to the orientation of $\partial{D}$. Let $\ga^6=l_1\cup...\cup l_6$ so that, for $i=1,...,6$,
$$\partial{D}\cap l_i=\{t_0\}\times\{p_i\}.$$

Assume that the annular neighborhood $A(\ga)$ of $\ga\subset \partial{V}=S^1\times \partial{D}^2$ is of the form
$$A(\ga)=\mathop{\bigcup}_{i=1}^6S^1\times[p_i-\varepsilon,p_i+\varepsilon],$$
for some small enough fixed constant $\varepsilon>0$. Let $T$ be an auxiliary surface of $(M,\ga^6)$, which consists of three disjoint annuli:
$$T=A_1\cup A_2\cup A_3,$$
where, for $i=1,2,3$, $A_i$ has the form
$$A_i=S^1_i\times [-1,1].$$
Choose an orientation reversing diffeomorphism
$f:\partial{T}\ra \ga$ so that
$$f(S^1_1\times\{1\})=l_1,~f(S^1_1\times\{-1\})=l_2,~f(S^1_2\times\{1\})=l_3,$$
$$f(S^1_2\times\{-1\})=l_6,~f(S^1_3\times\{1\})=l_4,~{\rm and~}f(S^1_3\times\{-1\})=l_5.$$

\begin{figure}[h]
\centering
\begin{overpic}[width=5.0in]{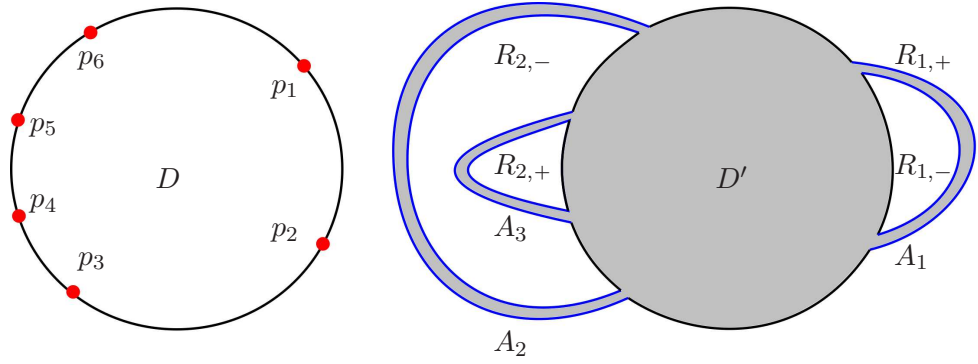}
    \put(15,15){$D$}
    \put(73,15){$D'$}
    \put(27,25){$p_1$}
    \put(27,10){$p_2$}
    \put(7,7){$p_3$}
    \put(7,28){$p_6$}
    \put(2,13){$p_4$}
    \put(2,21){$p_5$}
    \put(50,16.5){$R_{2,+}$}
    \put(50,28){$R_{2,-}$}
    \put(91.5,16.5){$R_{1,-}$}
    \put(91.5,28){$R_{1,+}$}
    \put(50,10){$A_3$}
    \put(50,-2){$A_2$}
    \put(91.5,7){$A_1$}
\end{overpic}
\vspace{0.05in}
\caption{A cross section $\{t\}\times D^2$ of the solid torus$V=S^1\times D^2$. The (red) dots in the left sub-figure represent the suture and the strips (with blue boundary) in the right sub-figure represent the three annuli $A_1,A_2$ and $A_3$. The shaded region is precisely the surface $D'$.}\label{fig_attach_annulus}
\end{figure}

Let
$$\widetilde{V}=V\mathop{\cup}_{f\times id} T\times[-\varepsilon,\varepsilon],$$
then $\widetilde{V}$ has four boundary components:
$$\widetilde{V}=R_{1,+}\cup R_{2,+}\cup R_{1,-}\cup R_{2,-},$$
so that $S^1_1\times\{\pm\varepsilon\}\subset R_{1,\pm}$ and $S^1_3\times\{\pm\varepsilon\}\subset R_{2,\pm}$. Suppose, for $i=1,2,3$, $S^1_i$ has the coordinate $t^i$, and $t^i_0$ is identified with $t_0\in S^1$ by $f$. Let
$$D'=D\cup (\{t^1_0\}\cup\{t^2_0\}\cup\{t^3_0\})\times[-1,1]\times[-\varepsilon,\varepsilon].$$
Thus, for $j=1,2$, we have $D'\cap R_{j,\pm}=C_{j,\pm}$. See Figure \ref{fig_attach_annulus}. Choose an orientation preserving diffeomorphism
$$h:(R_{1,+}\sqcup R_{2,+})\ra R_{1,-}\sqcup R_{2,-}$$
so that, for $j=1,2$,
$$h(C_{j,+})=C_{j,-}.$$
We can use $h$ to close up $\widetilde{V}$  to obtain a closure $(Y^{(6)},R^{(6)})$ of $(V,\ga^6)$. The surface $D'$ becomes closed oriented surface $\bar{D}^{(6)}$ of genus 2 inside $Y$. Define
$${SHM}(V,\ga^6,i)=\mathop{\bigoplus}_{\substack{\mathfrak{s}\in\mathfrak{S}(Y|R),\\ c_1(\mathfrak{s})[\bar{D}^{(6)}]=2i.}}\widecheck{HM}(Y,\mathfrak{s};\mathbb{Q}).$$ 
We know that
$${SHM}(V,\ga^6)\cong\mathop{\oplus}_{i\in\intg}{SHM}(V,\ga^6,i).$$

If we decompose the balanced sutured manifold $(V,\ga)$ along $D$, then, the result is a $3$-ball with one simple closed curve as the suture on its boundary. Thus, by Proposition 6.9 in Kronheimer and Mrowka \cite{kronheimer2010knots}, we know that
$${SHM}(V,\ga^6,1)\cong \mathbb{Q}.$$

On the other hand, we can also decompose $(V,\ga)$ along $-D$. A similar argument shows that
$${SHM}(V,\ga^6,-1)\cong\mathbb{Q}.$$
Hence, with $\mathbb{Q}$ coefficients, the rank of ${SHM}(V,\ga^6)$ is at least $2$. From Formula (\ref{eq_inductive_formula}) and the universal coefficient theorem, we know that ${SHM}(V,\ga^4)$, with either 
$\mathbb{Q}$ or $\intg$ coefficients, has a rank of at least two.

\begin{figure}[h]
\centering
\begin{overpic}[width=5in]{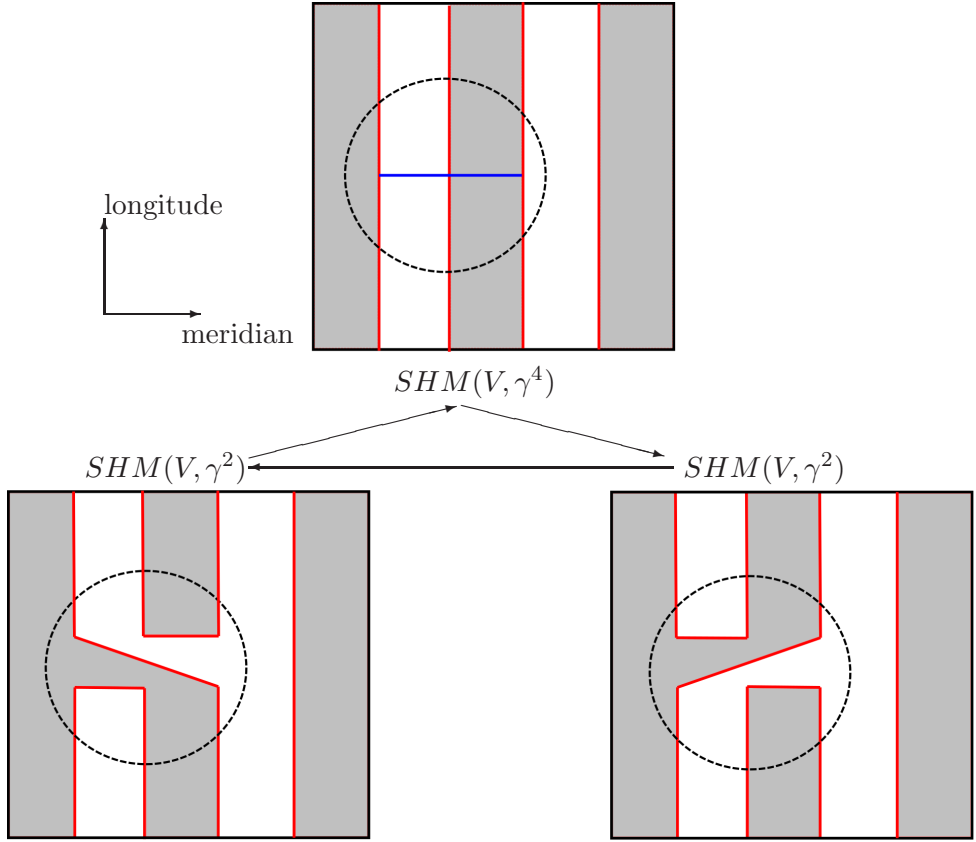}
    \put(10,55){\vector(1,0){10}}
    \put(10,55){\vector(0,1){10}}
    \put(18,52){meridian}
    \put(10,65.5){longitude}
    \put(40,47){${SHM}(V,\ga^4)$}
    \put(70,38){${SHM}(V,\ga^2)$}
    \put(8,38){${SHM}(V,\ga^2)$}
    \put(47,45.5){\vector(4,-1){21}}
    \put(69,39){\vector(-1,0){44}}
    \put(25,40){\vector(4,1){21.5}}
\end{overpic}
\vspace{0.05in}
\caption{The by-pass attachment along the horizontal (blue) arc $\al$. The change of sutures is limited in the dotted circles. The shaded region represents $R_-(\ga)$.}\label{fig_by_pass}
\end{figure}

To obtain an upper bound, we need to work with $\intg_2$ coefficients and use the by-pass exact triangles in sutured monopole Floer theory, which were introduced by Baldwin and Sivek in \cite{baldwin2016contact}. The by-pass attachments, as depicted in Figure \ref{fig_by_pass}, induce an exact triangle
\begin{equation*}\label{eq_exact_traingle_422}
\xymatrix{
&{SHM}(V,\ga^4)\ar[rd]^{\phi}&\\
{SHM}(V,\ga^2)\ar[ru]^{\rho}&&{SHM}(V,\ga^2)\ar[ll]^{\psi}
}
\end{equation*}
It is a basic fact that ${SHM}(V,\ga^2)\cong\intg_2$. So, with $\intg_2$ coefficients, the rank of ${SHM}(V,\ga^4)$ is at most two. By the universal coefficient theorem, we conclude that the rank with integral coefficients is also at most two.

It is also clear that, with $\intg$ coefficients, the rank of $SHM(V,\ga^4)$ is exactly two. Thus, by universal coefficients theorem, there is no even torsion in $G_{tor}$.
\epf

In \cite{kronheimer2010instanton}, Kronheimer and Mrowka constructed a particular closure of the manifold $(V,\ga^4)\sqcup (V,\ga^4)$. One can try to compute the monopole Floer homology of that closure directly, and we expect the following conjecture.

\begin{conj}
The torsion group $G_{tor}$ in Lemma \ref{lem_four_longitudes_Z_coefficients} is trivial.    
\end{conj}

\brem
In the proof of Lemma \ref{lem_four_longitudes_Z_coefficients}, we deal with $(V,\ga^6)$ instead of $(V,\ga^4)$ for the following reason: when dealing with $(V,\ga^4)$, we cannot pick a meridian disk $D$, which intersects with $\ga^4$ four times, and apply the construction, as in the proof of Lemma \ref{lem_four_longitudes_Z_coefficients}, to obtain a closed surface $\bar{D}$ in any closure $Y$ of $(V,\ga^4)$. This is due to a simple calculation of Euler characteristics. There is another subtlety in the construction of grading in the proof of Lemma \ref{lem_four_longitudes_Z_coefficients}. When pairing intersection points $p_1,...,p_6$, we did not simply pair the adjacent points, but, instead, we do it in a special way. This is because we want to have an even number of boundary components of the surface $D'$ on the positive or negative part of the boundary of $\widetilde{V}$.
\erem

\bthm\label{thm_solid_torus_with_longitude_sutures}
Suppose that $n=2k+1$ is odd. Then, with $\mathbb{Q}$ coefficients, there is a grading on ${SHM}(V,\ga^{2n})$ induced by a meridian disk of $V$ so that, with respect to this grading, we have
$${SHM}(V,\ga^{2n},i)\cong H_{i+k}(T^{n-1};\mathbb{Q})$$
for $-k\leq i\leq k$, and ${SHM}(V,\ga^{2n},i)=0$ for $|i|>k$. Here $T^{n-1}$ is the $(n-1)$-dimensional torus.
\ethm

\bpf
The basic case is trivial: if $k=0$, we have
$${SHM}(V,\ga^2)={SHM}(V,\ga^2,0)\cong \mathbb{Q}\cong H_0(T^0=\{pt\};\mathbb{Q}).$$

When $k=1$, the grading is already constructed in the proof of lemma \ref{lem_four_longitudes_Z_coefficients}, and we have
$${SHM}(V,\ga^6,\pm1)\cong\mathbb{Q}\cong H_0(T^2;\mathbb{Q})\cong H_2(T_2;\mathbb{Q}).$$

From the adjunction inequality (see Subsection 2.4 in \cite{kronheimer2010knots}), we know that, for $|i|>1$,
$${SHM}(V,\ga^6,i)=0.$$
Via Lemma \ref{lem_four_longitudes_Z_coefficients} and Formula (\ref{eq_inductive_formula}), we know that ${SHM}(V,\ga^6)\cong\mathbb{Q}^4$. Thus, we have
$${SHM}(V,\ga^6,0)\cong \mathbb{Q}^2\cong H_1(T^2;\mathbb{Q}).$$

For a general $k$, we argue in a similar way as we did for $(V,\ga^6)$ in the proof of Lemma \ref{lem_four_longitudes_Z_coefficients}. Let $D=\{t_0\}\times D^2$ be a meridian disk of $V$, and let
$$\partial{D}\cap \ga=\{p_1,...,p_{2n}\}.$$
The points $P_i$ are indexed according to the orientation of $\partial{D}$. The suture $\ga^{2n}$ can now be described as
$$\ga^{2n}=\bigcup_{i=1}^{2n}S^1\times\{p_i\}.$$
Pick an auxiliary surface $T$ for $(V,\ga^{2n})$, which consists of $n$ disjoint annuli:
$$T=\bigcup_{i=1}^n A_i.$$
Choose an orientation reversing diffeomorphism $f:\partial{T}\ra \ga$ so that 
$$f(\partial{A}_1)=S^1\times\{p_1,p_2\},$$
and, for $j=1,...,k$, we have
$$f(\partial{A}_{2j})=S^1\times\{p_{4k-1},p_{4k+2}\},~{\rm and}~f(\partial A_{2j+1})=S^1\times\{p_{4k},p_{4k+1}\}.$$
Let
$$\widetilde{V}=V\mathop{\cup}_{f\times id}T\times[-\varepsilon,\varepsilon],$$
we know that
$$\partial{\widetilde{V}}=\bigcup_{i=1}^{k+1}(R_{i,+}\cup R_{i,-}),$$
so that, for $j=1,...,k+1$,
$$A_{2j-1}\times\{\pm\varepsilon\}\subset R_{j,\pm}.$$
The meridian disk $D$ extends to a surface $D'\subset \widetilde{V}$ so that, for $j=1,...,k+1$,
$$\partial{D'}\cap R_{j,\pm}=C_{j,\pm}.$$
Choose an orientation preserving diffeomorphism
$$h:(R_{1,+}\cup...\cup R_{k+1,+})\ra R_{1,+}\cup...\cup R_{k+1,+}$$
so that, for $j=1,...,k+1$,
$$h(C_{j,+})=C_{j,-}.$$
Then, we get a closure $(Y^{(2n)},R^{(2n)})$ for $(V,\ga^{2n})$, and $D'$ becomes a closed oriented surface $\bar{D}^{(2n)}\subset Y^{(2n)}$.

Thus, we define a grading on ${SHM}(V,\ga^{2n})$ as follows:
$${SHM}(V,\ga^{2n},i)=\mathop{\bigoplus}_{\substack{\mathfrak{s}\in\mathfrak{S}(Y^{(2n)}|R^{(2n)}),\\ c_1(\mathfrak{s})[\bar{D}^{(2n)}]=2i.}}\widecheck{HM}(Y^{(2n)},\mathfrak{s};\mathbb{Q}).$$

Note $D'$ is obtained from $D$ by attaching $2k+1$ strips, so
$$\chi(\bar{D}^{2n})=\chi(D')=\chi(D)-(2k+1)=-2k.$$
By adjunction inequality, we know that
$${SHM}(V,\ga^{2n},i)=0$$
if $|i|>k$.

To compute the homology for each grading, we need to apply Floer excisions. Let $q_1,q_2\in \partial{D}\cap C_{1,+}\subset\partial{D}'$ be a pair of points. Suppose $q_1'=h^{-1}(q_1)$ and $q'_2=h^{-1}(q_2)$ where $h$ is the diffeomorphism we use to obtain the closure $(Y^{(2n)},R^{(2n)})$ for $(V,\ga^{2n})$. Suppose further we choose such an $h$ so that the following is true. 

(1) We have $q'_1,q'_2\in  \partial{D}\cap C_{1,+}\subset D'$.

(2) We have that $q'_1$ lies in between $p_6$ and $p_7$ and $q_2'$ lies in between $p_2$ and $p_3$.

(3) For $i=1,2$, we have
$$h(S^1\times\{q_i'\})= S^1\times\{q_i\}.$$
The three conditions above can be achieved by an $S^1$-invariant $h$. Pick two arcs $\be_1,\be_2\subset D$ so that, for $i=1,2$,
$$\be_i\cap\partial{D}=\partial{\be_i}=\{q_i,q'_i\}.$$
In $Y^{(2n)}$, $\be_1$ and $\be_2$ become two closed curves $\tilde{\be}_1$ and $\tilde{\be}_2$, respectively. Thus, there are two tori $T_1=\tilde{\be}_1\times S^1$ and $T_2=\tilde{\be}_2\times S^1$ inside $Y^{(2n)}$. Pick a $1$-cycle $\eta\subset R^{(2n)}\subset Y^{2n}$ to be the union of all images of $C_{i,+}\subset\partial{\widetilde{V}}$ in $Y^{(2n)}$. Clearly $\eta$ intersects both $T_1$ and $T_2$ transversely once.

\begin{figure}[h]
\centering
\begin{overpic}[width=5.0in]{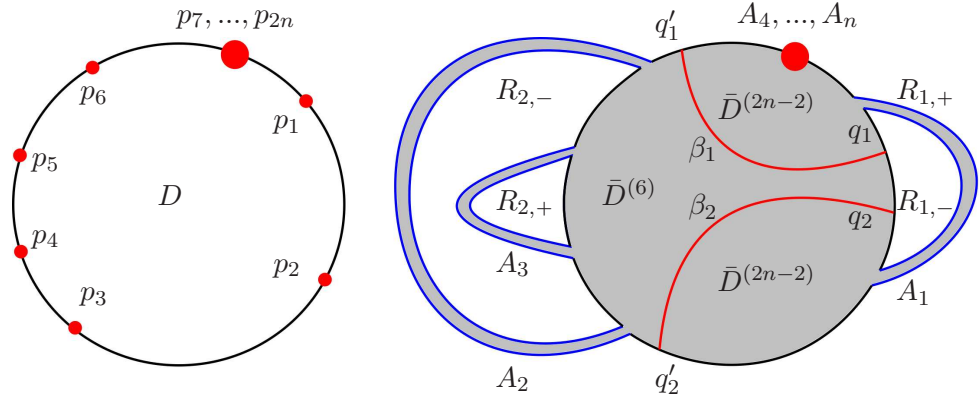}
    \put(15,17){$D$}
    \put(61,17){$\bar{D}^{(6)}$}
    \put(73,26){$\bar{D}^{(2n-2)}$}
    \put(73,8){$\bar{D}^{(2n-2)}$}
    \put(86.5,24){$q_1$}
    \put(86.5,15){$q_2$}
    \put(66.5,35){$q_1'$}
    \put(66.5,-2){$q_2'$}
    \put(27,25){$p_1$}
    \put(27,10){$p_2$}
    \put(7,7){$p_3$}
    \put(7,28){$p_6$}
    \put(2,13){$p_4$}
    \put(2,21){$p_5$}
    \put(50,16.5){$R_{2,+}$}
    \put(50,28){$R_{2,-}$}
    \put(91.5,16.5){$R_{1,-}$}
    \put(91.5,28){$R_{1,+}$}
    \put(50,10){$A_3$}
    \put(50,-2){$A_2$}
    \put(91.5,7){$A_1$}
    \put(17,36){$p_{7},...,p_{2n}$}
    \put(75,36){$A_4,...,A_{n}$}
    \put(70,22){$\be_1$}
    \put(70,16){$\be_2$}
\end{overpic}
\vspace{0.05in}
\caption{A cross section $\{t\}\times D^2$ of the solid torus$V=S^1\times D^2$. The (red) dots in the left sub-figure represent the suture and the stripes (with blue boundary) in the right sub-figure representing the three annuli $A_1,A_2$ and $A_3$. $A_4,...,A_{n}$ are not depicted. The two (red) arcs inside $D$ are $\be_1,\be_2$. The shaded region represents the surface $D'$.}\label{fig_cut_along_annulus2}
\end{figure}

We can perform a Floer excision along $T_1$ and $T_2$, or to be more precise, the inverse operation of the Floer excision introduced in Subsection \ref{subsec_sutured_monopoles}. The result of this reversed Floer excision is a disjoint union of two $3$-manifolds, $Y^{(2n-4)}$ and $Y^{(6)}$. Along this process, the surface $R^{(2n)}$ is cut into $R^{(2n-4)}\cup R^{(6)}$, and the surface $\bar{D}^{(2n)}$ is cut into $\bar{D}^{(2n-4)}\cup\bar{D}^{(6)}$. See Figure \ref{fig_cut_along_annulus2}.

\vspace{0.2in}

As described in Subsection \ref{subsec_sutured_monopoles}, there is a cobordism $W$ from $Y^{(2n-4)}\sqcup Y^{(6)}$ to $Y^{(2n)}$. Inside the cobordism $W$, there is a $3$-dimensional cobordism between $\bar{D}^{(2n-4)}\sqcup \bar{D}^{(6)}\subset Y^{(2n-4)}\sqcup Y^{(6)}$ and $\bar{D}^{(2n)}\subset Y^{(2n)}$. This $3$-dimensional cobordism is obtained from a $2$-dimensional analogue of the $3$-dimensional Floer excision. Thus, if $\mathfrak{s}$ is a spin${}^c$ structure on $W$ so that
$$c_1(\mathfrak{s})[\bar{D}^{(2n-4)}]=x,~c_1(\mathfrak{s})[\bar{D}^{(6)}]=y,$$
then we must have
$$c_1(\mathfrak{s})[\bar{D}^{(2n)}]=x+y.$$
As a result, there is a product formula that enables us to compute ${SHM}(V,\ga^{2n})$ out of ${SHM}(V,\ga^{2n-4})$ and ${SHM}(V,\ga^{6})$. After a grading shifting, this product formula is precisely the one we compute $H_{*}(T^{n-1})$ from $T^{n-1}=T^{n-3}\times T^2$. Thus, we conclude the proof of Theorem \ref{thm_solid_torus_with_longitude_sutures}.
\epf

One question arises in the proof of Theorem \ref{thm_solid_torus_with_longitude_sutures}. Fix a suitable field $\mathcal{F}$ of characteristic 2. There is a by-pass exact triangle for a general $(V,\ga^{2n})$:
\begin{equation*}\label{eq_exact_traingle_general}
\xymatrix{
&{SHM}(V,\ga^{2n})\ar[rd]^{\phi}&\\
{SHM}(V,\ga^{2n-2})\ar[ru]^{\rho}&&{SHM}(V,\ga^{2n-2})\ar[ll]^{\psi}
}
\end{equation*}
From Formula (\ref{eq_inductive_formula}), we know that, for $n>1$,
$${SHM}(V,\ga^{2n})\cong (\mathcal{F})^{2^{n-2}.}$$
This forces the map $\psi$ to be $0$. Hence, $\rho$ is injective and $\phi$ is surjective. If we assume that $n=2$, then we know from Geiges \cite{geiges2008introduction} that there is a unique tight contact structure $\xi_0$ that is compatible with $(V,\ga^2)$. From \cite{baldwin2016contact}, we know that the contact element of $\xi_0$ generates ${SHM}(V,\ga^2)\cong\mathcal{F}$. Since the map associated to the by-pass attachment preserves contact elements, we know that after attaching the by-pass associated to $\psi$, $\xi_0$ becomes overtwisted, and after attaching the by-pass associated to $\rho$, $\xi_0$ becomes a compatible contact structure $\xi_1$ on $(V,\ga^4)$ so that the contact element of $\xi_1$ generates $im(\rho)\cong\mathcal{F}\subset {SHM}(V,\ga^4)$. If there was another compatible contact structure $\xi_2$ on $(V,\ga^4)$ so that, after the by-pass associated to $\phi$, it became $\xi_0$ on $(V,\ga^2)$, then we could conclude that ${SHM}(V,\ga^4)$ was simply generated by the two contact elements of $\xi_1$ and $\xi_2$. Then, we could also look at a general $(V,\ga^{2n})$. However, by-pass attachments do not necessarily have inverses, and this lead to the following question:

\begin{quest}
Is ${SHM}(V,\ga^{2n})$ fully generated by contact elements of compatible contact structures?
\end{quest}

\subsection{The connected sum formula}\label{subsec_connected_sum}
In this subsection, we derive the connected sum formula for sutured monopole Floer homology. First, we prove the following proposition.

\bprop\label{prop_oriented_skein_relation_for_shm}
We use $\intg_2$ coefficients. Suppose three oriented links $L_0$, $L_1$ and $L_2$ are the same outside a $3$-ball $B^3$, and, inside $B^3$, they are depicted as in Figure \ref{fig_oriented_skein_relation}. We have the following. 

(1) If $L_2$ has one more component than $L_0$ and $L_1$, then there is an exact triangle:
\begin{equation*}
\xymatrix{
KHM(S^3,L_0)\ar[rr]&&KHM(S^3,L_1)\ar[ld]\\
&KHM(S^3,L_2)\ar[lu]&
}    
\end{equation*}

(2) If $L_2$ has one less component than $L_0$ and $L_1$, then there is an exact triangle:

\begin{equation*}
\xymatrix{
KHM(S^3,L_0)\ar[rr]&&KHM(S^3,L_1)\ar[ld]\\
&KHM(S^3,L_2)\otimes(\intg_2)^{4}\ar[lu]&
}    
\end{equation*}
\eprop

\bpf
It follows from an analogous argument in sutured instanton Floer theory in Kronheimer and Mrowka \cite{kronheimer2010instanton}. We sketch the proof as follows: The monopole knot Floer homology $KHM(S^3,L_0)$ for a link $L_0\subset S^3$ is defined by taking the sutured monopole Floer homology of the balanced sutured manifold $(S^3(L_0),\Gamma_{\mu})$, where $S^3(L_0)$ is the link complement, and $\Gamma_{\mu}$ consists of a pair of meridians on each boundary component of $S^3(L_0)$. Let $(Y_0,R)$ be a closure of $(S^3(L_1),\Gamma_{\mu})$. To obtain an exact triangle for the oriented skein relation, we pick a small circle $\al$ linking around the crossing, on which we perform the crossing change and oriented smoothing. The curve $\al$ naturally embeds into $Y_0$ and is disjoint from the surface $R\subset Y_0$. Then, from Kronheimer, Mrowka, Ozsv\'ath, and Szab\'o \cite{kronheimer2007monopolesandlens}, there is a surgery exact triangle relating the monopole Floer homologies of the $3$-manifolds obtained by performing $(-1)$, $0$, and $\infty$ surgery along the curve $\al\subset Y_0$. There is a canonical framing for the curve $\al\subset S^3$, and the surgeries slopes are the ones with respect to the canonical framing. 

When performing the $\infty$ surgery, we obtain $Y_0$. When performing the $(-1)$ surgery, we obtain a closure of $(S^3(L_2),\Gamma_{\mu})$, which gives rise to $KHM(S^3,L_2)$. When performing the $0$ surgery, we get a closure of the balanced sutured manifold $(M,\ga)$, which is obtained from $(S^3(L_0),\Gamma_{\mu})$ by performing a $0$-surgery along $\al$. To further relate $(M,\ga)$ to $L_0$, we need to perform a sutured manifold decomposition of $(M,\ga)$, along an annulus $A$ arising from $\al$ and the $0$-surgery, as explained in Kronheimer and Mrowka \cite{kronheimer2010instanton}. Suppose that $(M',\ga')$ is obtained from $(M,\ga)$ by such a decomposition along $A$. There are two cases.

{\bf Case 1.} When $L_2$ has one more component than $L_0$ and $L_1$, then $(M',\ga')$ is diffeomorphic to $(S^3(L_2),\Gamma_{\mu})$. Hence, we are done.

{\bf Case 2.} When $L_2$ has one less component than $L_0$ and $L_1$, then $(M',\ga')$ is $(S^3(L_2),\Gamma_{\mu})$ except that on one component of $\partial S^3(L_2)$, there are six meridians as the suture rather than two. So, to obtain $(S^3(L_2),\Gamma_{\mu})$, we need to further decompose $(M',\ga')$ by another annulus, just as we did in the proof of Lemma \ref{lem_four_longitudes_Z_coefficients}, where we obtained two copies of $(V,\ga^4)$ from $(V,\ga^6)$. The result of this second sutured manifold decomposition is a disjoint union of $(S^3(L_2),\Gamma_{\mu})$ with $(V,\ga^6)$. Hence, we are done.
\epf

As a corollary to Proposition \ref{prop_oriented_skein_relation_for_shm}, we derive the following corollary, independent of the works by \cite{fintushel1998knots} or \cite{meng1996underline}.

\bcor
With $\intg_2$ coefficients and the canonical $\intg_2$ grading of monopole Floer homology, the Euler characteristics of ${\rm KHM}(S^3,K,i)$ (for definition, see \cite{kronheimer2010knots}) corresponds to the coefficients of a suitable version of Alexander polynomial of the knot $K\subset S^3$.
\ecor

\bpf
It follows from an analogous argument in sutured instanton Floer theory in \cite{kronheimer2010instanton}.
\epf

Suppose that $Y$ is a closed oriented $3$-manifold. Let $Y(n)$ denote the manifold obtained by removing $n$ disjoint $3$-balls from $Y$. We can make $Y(n)$ to be a balanced sutured manifold $(Y(n),\delta^n)$, where $\delta^n$ consists of one simple closed curve on each boundary sphere of $Y(n)$.

The following two lemmas are straightforward. 
\blem\label{lem_topological_description_for_digging_one_more_ball}
Suppose $Y$ is a closed oriented $3$-manifold and $n\in \intg$ is no less than $2$, then 
$$Y(n)\cong (Y(n-1)\sqcup S^3(2),\delta^{n-1}\cup\delta^2)\cup h,$$
where $h=(\phi,S,D^3,\delta)$ is a contact $1$-handle so that $\phi$ sends one component of $S$ to $\partial{Y(n-1)}$ and the other component to $\partial{S^3(2)}.$
\elem

\blem\label{lem_topological_description_for_connected_sum}
Suppose $(M_1,\ga_1)$ and $(M_2,\ga_2)$ are two balanced sutured manifolds. Also, suppose that $(S^3(2),\delta^2)$ is defined as above, and its two boundary components are
$$\partial{S^3(2)}=S^2_1\cup S^2_2.$$
Then, we have
$$(M_1\sharp M_2,\ga\cup\ga_2)\cong(M_1\sqcup M_2\sqcup S^3(2),\ga_1\cup \ga_2\cup\delta^2)\cup h_1\cup h_2.$$
Here, for $i=1,2,$ $h_i=(\phi_i,S_i,D^3_i,\delta_i)$ is a contact $1$-handle so that $\phi_i$ maps one component of $S_i$ to $\partial{M}_i$ and the other component of $S_i$ to $S_i^2$.
\elem

\brem
In Lemmas \ref{lem_topological_description_for_digging_one_more_ball} and \ref{lem_topological_description_for_connected_sum}, we do not require a sutured manifold $(M,\ga)$ to have a global compatible contact structure. However, we can identify a collar of the boundary of $M$ with $\partial{M}\times [0,1]$, and assume that there is an $I$-invariant contact structure in the collar so that, under this contact structure, $\partial{M}$ is a convex surface with $\ga$ being the dividing set. Thus, the contact handle attachment makes sense.
\erem

From Lemmas \ref{lem_topological_description_for_digging_one_more_ball} and \ref{lem_topological_description_for_connected_sum}, we can see the significant role played by $(S^3(2),\delta^2)$. So, we proceed to compute its sutured monopole Floer homology. 

\blem\label{lem_injective_maps_exists}
For any closed $3$-manifold $Y$ and every positive integer $n$, there is an injective map
$${SHM}(Y(n),\delta^n)\ra {SHM}(Y(n+1),\delta^{n+1}).$$
\elem
\bpf
We can obtain $(Y(n+1),\delta^{n+1})$ from $(Y(n),\delta^n)$ by attaching a contact $2$-handle. If we further attach a contact $3$-handle to it, the result will be $(Y(n),\delta^n)$ again. The pair of handles forms a 2-3 cancelation pair, as in Li \cite{li2018gluing}. Thus, the composition is the identity, and the $2$-handle attachment induces the desired injective map.
\epf

\bcor\label{cor_SHM_of_S_three_two}
We have ${SHM}(S^3(2),\delta^2;\intg_2)\cong (\intg_2)^2.$
\ecor

\bpf
It follows from the proof of an analogous statement in sutured instanton Floer theory in Baldwin and Sivek \cite{baldwin2016instanton}. A sketch of the proof is as follows: First, as in \cite{baldwin2016contact}, there is an isomorphism
$$SHM(S^3(2),\delta^2)\cong KHM(S^3,U_2),$$
where $U_2$ is the unlink with two components. Note that there is an oriented skein relation, which relates two copies of the unknot $U_1$ and one copy of $U_2$. Thus, from Proposition \ref{prop_oriented_skein_relation_for_shm} we have
\begin{equation*}
    \xymatrix{
    KHM(S^3,U_1)\ar[rr]&&KHM(S^3,U_1)\ar[ld]\\
    &KHM(S^3,U_2)\ar[lu]&
    }
\end{equation*}
Since $KHM(S^3,U_1)\cong\intg_2$, we know that
$$KHM(S^3,U_2)\cong(\intg_2)^2~{\rm or}~0.$$
The second possibility is then ruled out by Lemma \ref{lem_injective_maps_exists} since
$$SHM(S^3(1),\delta^1)\cong\intg_2.$$
So, we conclude that
$$SHM(S^3(2),\delta^2)\cong KHM(S^3,U_2)\cong(\intg_2)^2.$$
\epf

\bcor\label{cor_connected_sum_formula_for_SHM}
Suppose $(M_1,\ga_1)$ and $(M_2,\ga_2)$ are balanced sutured manifolds. Then, we have
$${SHM}(M_1\sharp M_2,\ga_1\cup\ga_2)\cong {SHM}(M_1\sqcup M_2,\ga_1\cup\ga_2)\otimes (\intg_2)^2.$$
\ecor

\bpf
This follows directly from Lemma \ref{lem_topological_description_for_connected_sum} and Corollary \ref{cor_SHM_of_S_three_two}.
\epf

\bcor
Suppose $L$ is a link in $S^3$. Then, with any coefficients, ${\rm KHM}(S^3,L)\neq0$.
\ecor

\bpf
Lemma \ref{lem_injective_maps_exists} makes sure that $SHM(S^3(2),\delta^2)$ has a rank of at least one with any coefficients. If $L$ is non-splitting, then the balanced sutured manifold $(S^3(L),\Gamma_{\mu})$ is taut, and the non-vanishing statement follows from Kronheimer and Mrowka \cite{kronheimer2010knots}. If $L$ has separable components, we can apply Lemma \ref{lem_topological_description_for_connected_sum}.
\epf

The discussion in the instanton settings would be completely analogous. We use the field of complex numbers $\mathbb{C}$ as coefficients and have the following proposition.
\bprop
Suppose $(M_1,\ga_1)$ and $(M_2,\ga_2)$ are two balanced sutured manifolds, then
$${SHI}(M_1\sharp M_2,\ga_1\cup\ga_2)\cong { SHI}(M_1,\ga_1)\otimes { SHI}(M_2,\ga_2)\otimes \mathbb{C}^2.$$
\eprop

This formula can also be applied to the framed instanton Floer homology of closed $3$-manifolds. Suppose $Y$ is a closed oriented $3$-manifold, then we can connected sum $Y$ with $T^3$, and pick $\omega$ to be a circle that represents a generator of $H_1(T^3)$. The pair $(Y\sharp T^3,\omega)$ is then admissible and we can form the {\it framed instanton Floer homology} of $Y$:
$$I^{\sharp}(Y)=I^{\omega}(Y\sharp T^3).$$
In \cite{kronheimer2010knots}, Kronheimer and Mrowka discuss the relation between the framed instanton Floer homology of a closed $3$-manifold and the sutured instanton Floer homology of $(Y(1),\delta^1)$. As a corollary to the connected sum formula for sutured instanton Floer theory, we have the following: 

\bcor
Suppose $Y_1$ and $Y_2$ are two closed oriented $3$-manifolds. Then, as vector spaces over complex numbers, we have
$$I^{\sharp}(Y_1)\otimes I^{\sharp}(Y_2)\cong I^{\sharp}(Y_1\sharp Y_2).$$
\ecor

%————End from here————

\bibliography{Index}%for reference
\end{document}